\documentclass[conference]{IEEEtran}
\IEEEoverridecommandlockouts
\usepackage{cite}
\usepackage{amsmath,amssymb,amsfonts,amsthm}
\usepackage{hyperref,ulem,comment}
\usepackage[ruled,linesnumbered]{algorithm2e}
\usepackage[switch]{lineno}
\usepackage{graphicx}
\usepackage{textcomp}
\usepackage{xcolor}
\def\BibTeX{{\rm B\kern-.05em{\sc i\kern-.025em b}\kern-.08em
    T\kern-.1667em\lower.7ex\hbox{E}\kern-.125emX}}
    
\usepackage{graphicx}
\usepackage{subfigure}
\usepackage{float}
\usepackage{multirow}
\usepackage{booktabs}
\usepackage{verbatim}
\usepackage{amsthm}

\theoremstyle{definition}
\newtheorem{example}{Example}[section]
\newtheorem{remark}{Remark}[section]

\newcommand{\TV}{\operatorname{TV}}
\newcommand{\NN}{\operatorname{NN}}

\newcommand{\abs}[1]{\left\vert#1\right\vert}

\newcommand{\wt}[1]{{\widetilde{#1}}}

\newcommand{\JR}{ \mathcal{J}_{\textrm{R}} }
\newcommand{\JG}{ \mathcal{J}_{\textrm{G}} }
    
\DeclareMathOperator{\Tr}{Tr}

\begin{document}

% \title{Roughness Index for Understanding Loss Landscapes of Neural Network Models 
% of    Partial  Differential Equations*\\
% {\footnotesize \textsuperscript{*}Note: Sub-titles are not captured in Xplore and
% should not be used}
% \thanks{Identify applicable funding agency here. If none, delete this.}
% }

\title{Roughness Index for  Loss Landscapes of Neural Network Models of Partial  Differential Equations
}

\author{\IEEEauthorblockN{Keke Wu}
\IEEEauthorblockA{\textit{School of Mathematical Sciences} \\
\textit{Shanghai Jiao Tong University}\\
Shanghai, China \\
wukekever@sjtu.edu.cn}
\and
\IEEEauthorblockN{Xiangru Jian}
\IEEEauthorblockA{\textit{School of Data Science} \\
\textit{City University of Hong Kong}\\
Hong Kong, China \\
xiangjian2-c@my.cityu.edu.hk}
\and
\IEEEauthorblockN{Rui Du}
\IEEEauthorblockA{\textit{School of Mathematical Sciences} \\
\textit{Soochow University}\\
Suzhou, China \\
durui@suda.edu.cn}
\and
\IEEEauthorblockN{Jingrun Chen}
\IEEEauthorblockA{\textit{School of Mathematical Sciences and Suzhou Institute for Advanced Research} \\
\textit{University of Science and Technology of China}\\
Hefei, China \\
jingrunchen@ustc.edu.cn}
\and
\IEEEauthorblockN{  Xiang ZHOU}
\IEEEauthorblockA{\textit{School of Data Science and Department of Mathematics } \\
\textit{ City University of Hong Kong}\\
Hong Kong, China \\
xizhou@cityu.edu.hk}
 
}

\maketitle

\begin{abstract}
Loss landscape is a useful tool to characterize and compare neural network models. The main challenge for analysis of loss landscape for the deep neural networks is that they are generally highly non-convex in very high dimensional space. In this paper, we develop the “roughness” concept for understanding such landscapes in high dimensions and apply this technique to study two neural network models arising from solving differential equations. Our main innovation is the proposal of a well-defined and easy-to-compute roughness index (RI) which is based on the mean and variance of the (normalized) total variation for one-dimensional functions projected on randomly sampled directions. A large RI at the local minimizer hints an oscillatory landscape profile and indicates a severe challenge for the first-order optimization method. Particularly, we observe the increasing-then-decreasing pattern for RI along the gradient descent path in most models. We apply our method to two types of loss functions used to solve partial differential equations (PDEs) when the solution of PDE is parametrized by neural networks. Our empirical results on these PDE problems reveal important and consistent observations that the landscapes from the deep Galerkin method around its local minimizers are less rough than the deep Ritz method.\end{abstract}


\section{Introduction}

In recent years, solving partial differential equations by deep neural networks (DNNs) has brought significant interests from the community of scientific computing; see \cite{beck2020overview} for reviews and references therein. Due to its powerful representation ability, a DNN can well approximate a target function in high dimensions. Given a PDE, the basic idea is to use a DNN as the trial function to approximate the PDE solution. The optimal set of parameters in the DNN is obtained by minimizing a loss function in different forms\cite{E2018, deepGalerkin2018, raissi2019physics}. Since the loss function lives in the high-dimensional parameter space and is highly nonconvex, it is difficult to find the global minimizer. The minimization problem is often solved by the stochastic gradient descent method\cite{BottouSGD}. The complexity of loss landscapes makes the training process and the numerical results highly depend on the DNN structure, the optimization method as well as  the initialization\cite{choromanska2015open}.

Efforts towards to better understandings of loss landscapes include studies on specific problems \cite{nguyen2018loss,gamst2017energy}, geometry of local minima \cite{swirszcz2016local,pmlr-v70-dinh17b,Baldassi161}, energy barriers\cite{draxler2018essentially}, mean field limit \cite{mei2018mean}, as well as neural tangent kernel limit \cite{jacot2018neural}.
Due to the high dimensionality of the parameter space, it is  difficult to visualize the loss function. One strategy is to project the loss function onto a low-dimensional space with the random choice of directions and filter-wise normalization \cite{li2018visualizing}. This has been used to show the advantage of some residual NNs \cite{He2015} over fully connected NNs. In addition, 
the volume of basin of attractor has been considered to characterize the flatness of minima\cite{wu2017towards}.

Our interest is how to  understand and compare  two loss functions in the background  of solving PDEs.
In  this PDE context,  one have the same network architectures and the training data  to solve the same PDE, but 
have different  forms of the loss functions. Since both  loss functions solve  the same PDE,   we can fairly compare
the performance of two loss functions in this task. 
% In addition, the loss functions are  population risk and thus free of any particular dataset.
% Thus, we have an ideal experimental problem to study the loss functions themselves by excluding any potential  ambiguity.
This paper considers two representative methods for  solving PDE with DNN. One is the variation-based model -- deep Ritz method (DRM) \cite{E2018} and the other one is the residual-based model -- deep Galerkin method (DGM) \cite{deepGalerkin2018}. %It should be noted that the DGM uses the mean squared PDE residual as the loss function, which is in fact different from the Galerkin approach \cite{zang2020weak}  in numerical PDEs. 

It is well-known \cite{Goodfellow2016,pmlr-v48-hardt16,2016arXiv160904836S} that the loss function is complex due to non-convexity,
and has many oscillatory local minimizers  in  the valley 
of a   ``good'' minimizer by SGD. Such  good minimizers  are conjectured to  be wide and flat in geometry.
and  thus  have better  generalization ability.
To reach such good (local) minimizers, the training process relies on the noise injected by the stochastic optimization method
to climb over the small barriers so as to achieve a better accuracy and generalization error.
Therefore, the landscape is  essentially {\it rough}
and the training process is an exploration process of the rough landscape before eventually 
hitting the final solution. 

In this paper, we propose a quantitive index to describe this roughness concept
and use it  to 
measure at an any point how  rough the landscapes are for the different models 
in solving the PDEs.
This index will be used to  characterize  the accumulated effect of the small-scale  oscillatory wells
within   neighborhoods of   numerically obtained minimizers. 
We  call  it the {\it roughness index} (RI).
This index is associated with  each minimizer  which is found by the standard stochastic optimization approach.
But meanwhile, this quantity is delocalized in the sense that it does not rely on the eigenvalues on the minimizer 
and  it is beyond    the infinitesimal quadratic approximation. 
This index may depend on the size of the neighborhood, which is a box 
in our computation. Ideally this length scale should be the typical size  scale of the basin of attraction.  We practically compute the index for varied size  and identify 
the consistent result within a range of proper size.

By computing the RI for various local minimizers
of DGM and DRM applied to  the Poisson equation,  we find the consistent and distinctive differences: 
the DGM's minimizers have a smaller RI while 
the DRM's minimizers have a larger RI. 
We also track the RI along the training  trajectory and  
find for typical initialized parameters in the NN, 
the roughness index is small, and  the DRM's roughness index 
gradually increase when approaching the minimizers.

In a nutshell, by studying the roughness index in the space of high dimensional parameter space, 
we can reveal a few interesting and phenomenal understandings about the loss landscape in quantitative ways  which have not yet been explored.  
This roughness index is not restricted to the NN models for PDEs, but a potential tool for analyzing general machine-learning landscapes.

This paper is organized as follows. We first give an introduction of methods for solving PDEs by DNNs: 
the DRM and DGM. Section \ref{sec:methods} is our main part  
to define and compute  roughness index.
Section \ref{sec:result}  applies the RI  to different models, different neural networks, and different dimensions and different  PDEs. 
Conclusive remarks are drawn in Section \ref{sec:conclusion}.

\section{Related works}\label{sec:methods}
\subsection{Solving PDEs by deep neural networks} \label{ssec:nn}
% In machine-learning tasks, the mean-square error between the labeled data and the model by a NN is often used as the loss function. 
When using a NN to solve a given PDE, there are multiple choices to construct the loss function. If the PDE can be derived as the Euler-Lagrange equation of a variational problem, then this variational problem can be defined as the loss function; see DRM \cite{E2018} for example. In contrast, DGM \cite{deepGalerkin2018}  has the loss function as the mean-square error or the residual associated to the given PDE.  For completeness, we shall first review these two methods for the elliptic equation where the variational loss function exists. 

Consider the Poisson equation over a bounded domain $\Omega\subset\mathbb{R}^d$ 
\begin{equation}
\begin{cases}
-\Delta u(x) = f(x), & \;\text{in} \;\Omega,\\
u(x) = g(x), & \; \text{on} \; \partial\Omega, 
\end{cases}
\label{eq}
\end{equation}
where $f, g$ are given functions.
Denote $u(x;\theta)$ the approximate NN solution with the set of parameters  $\theta$. 
% The collection of all $u(x;\theta)$ is used as the trial space.
The network structure employed here is ResNet\cite{He2015} with several residual blocks or a fully-connected Net (FCNet)\cite{Goodfellow2016}.
% ResNet is built by stacking some residual blocks with each block containing one input, two weight layers, two activation functions, one shortcut connection, and one output.
%see the right column of Figure \ref{fig: fcnAndresnet} for illustration. 
Consider a ResNet with $N$ residual blocks. For the $i$-th block, let $L^{i}[x] \in \mathbb{R}^{w\times1}$ be the input, $W_1^{i}, W_2^{i} \in \mathbb{R}^{w \times w} $ and $b_1^{i}, b_2^{i} \in \mathbb{R}^{w\times1}$ be the weight matrices and bias vectors, $\sigma(\cdot)$ be the activation function, then the output $L^{i+1}[x]$ can be written as 
\begin{equation*}
 {L^{i+1}[x] = L^{i}[x] + \sigma (W_2^{i}\cdot\sigma (W_1^{i}\cdot L^{i}[x] + b_1^{i}) + b_2^{i}), }
\end{equation*}
{where $i = 0, \cdots, N-1.$}
The input and the output are $L^0(x) = W^0 \cdot x + b^0$ and $ L^{N+1} [x] = W^{N+1} \cdot L^N [x]  + b^{N+1}$ with $W^0 \in \mathbb{R}^{w \times d}, b^0 \in  \mathbb{R}^{w \times 1}$ and $W^{N+1} \in \mathbb{R}^{1 \times w}, b^{N+1} \in  \mathbb{R}$. 
For the $i$-layer of FCNet with $2N+2$ layers, let $L^{i}[x] \in \mathbb{R}^{w\times1}$ be the input, $W^{i} \in \mathbb{R}^{w \times w} $ and $b^{i} \in \mathbb{R}^{w\times1}$ be the weight matrices and bias vectors, then the output $L^{i+1}[x]$ can be written as 
\begin{equation*}
L^{i+1}[x] = \sigma (W^{i}\cdot L^{i}[x] + b^{i}), i = 0, 1, \cdots, 2N.
\end{equation*}
The   input and   output are $L^0(x) = W^0 \cdot x + b^0$ and $ L^{2N+2} [x] = W^{2N+2} \cdot L^{2N+1} [x]  + b^{2N+2}$ with $W^0 \in \mathbb{R}^{w \times d}, b^0 \in  \mathbb{R}^{w \times 1}$ and $W^{2N+2} \in \mathbb{R}^{1 \times w}, b^{2N+2} \in  \mathbb{R}$. The number of neurons in each hidden layer (neural width) is $w$. Therefore, the total number of parameters in ResNet or FCNet is $2Nw^2 + (d + 2N+2)w + 1$.
Since the Hessian information needs to be calculated, we use the $swish$ function ($x(1+e^{-x})^{-1}$) as the activation function in what follows.
Boundary condition can be enforced exactly by constructing a special neural network. DGM and DRM only differ by their loss functions in terms of $\theta$, which are
\begin{equation}\label{eqn:lossDGM}
\mathcal{J}_{\textrm{G}} (\theta) = \int_{\Omega} {|-\Delta u(x;\theta) - f(x)|}^2 \mathrm{d}x,
\end{equation}
and
\begin{equation}\label{eqn:lossDRM}
\mathcal{J}_{\textrm{R}} (\theta) = \int_{\Omega} \left(\frac{1}{2}|\nabla u(x;\theta)|^2 - f(x)u(x;\theta)\right) \mathrm{d}x.
\end{equation}
% From these two loss functions, it is easy to see that DGM aims to minimize the residual  in the least-squares sense \cite{deepGalerkin2018} and the DRM minimizes the energy $\mathcal{J}_{\textrm{R}} (\theta)$ \cite{E2018}.

A minimizer is obtained by Adam optimizer. 
Derivatives of $u(x;\theta)$ are calculated by the automatic differentiation. Monte Carlo method is applied to approximate the integrals in DGM and DRM by $N$ samples. 
% Various choices of $N$ are tested in  numerical examples so that the sampling error does not affect any quantities of our interests.  
In 1D, instead, the Simpson's rule is used for better accuracy. One \textit{epoch} refers to the period of processing  $N$ samples, i.e., one time step in Adam. We typically set the batch size $N=200,1000,10000$ in the 1D, 3D, 10D PDE, respectively.

\subsection{Eigenvalue-based index}
The loss landscape is complicated and typically there are many minima of interest. For example, for a simple NN, the minima of loss function
may  lie in a very flat basin~\cite{gamst2017energy}. To understand the loss landscapes of DGM and DRM, we
first consider  the concept of ``volume of basin of attractor'' proposed in \cite{wu2017towards}.
Their use of the (Lebsque) measure of the basin for each attractor is 
an appealing idea. However, it is almost impossible in reality  to find the exact basin and precisely measure its volume 
in high dimensional space. As a compromise, 
\cite{wu2017towards} in fact used the Hessian matrix 
at the {minimum point} to represent the ``volume of the basin'' of this minimum point. 
Precisely, for a given minimizer $\theta^*$, one can compute the Hessian $H$ of loss function with respect to $\theta$ and evaluate it at $\theta^*$.
Since the volume of the sublevel set of a {\it quadratic} form is proportional   to the product of   eigenvalues,  \cite{wu2017towards} used
the logarithm of the product of top-$k$ eigenvalues ($k$ is truncated to keep only significant nonzero eigenvalues) of $H(\theta^*)$ to approximate the inverse volume of basin of attractor
\begin{equation}\label{basin}
V(k) := \sum_{i = 1}^k \log_{10} (\lambda_i(H(\theta^*))).
\end{equation}
\eqref{basin} provides a quantitative characterization of the size of the basin around a minimizer for the local quadratic approximation of the landscape. 
A  small $V$ means a   ``flat'' valley near $\theta^*$ and is regarded 
to have a large volume of basin, which arguably is 
able to  generalize well\cite{li2018visualizing,wu2017towards}. 
We emphasize  that the index $V$ in \eqref{basin} only relies on the Hessian information at the minimizer, thus is essentially  a local quantity for characterizing the flatness and the assumption behind is that the landscape around $\theta^*$ is convex and smooth. 
However, the neighboring region for such assumptions to be valid could be very small in practice and it is hard to justify the applicability of  
this index to represent the real non-convex  behaviors around the local minimum points.

\subsection{Normalized total variation for 1D functions}
Total variation (TV) is a commonly used norm in applied mathematics for regularity of a function.
For instance,  TV has been used in image denoising  
as a penalty to suppress the spurious detail \cite{RUDIN1992259,selesnick2012total}. It is also adopted in the statistical  learning for the purpose of smoothing and regularization in fitting data.
It is one of natural candidates to describe the ``regularity'' or ``roughness'' of the signals. 
We propose to utilize the concept of TV to construct  roughness index. 

Recall that the TV of a continuous function $f$ from $[a, b]$ to $\mathbb{R}$ is given by
$$
\operatorname{TV}(f)=\sup \sum_{k=0}^{n-1}\left|f\left(x_{k+1}\right)-f\left(x_{k}\right)\right|
$$
where the sup is taken over all possible partitions, $a=x_{0}<\ldots<x_{n}=b$.
If $f$ is absolutely continuous, we can write
$$\operatorname{TV}(f)=\int_a^b |f'(x)|\mathrm{d}x.$$

\begin{figure}[htbp!]
	\centering
	\includegraphics[width=0.5\textwidth]{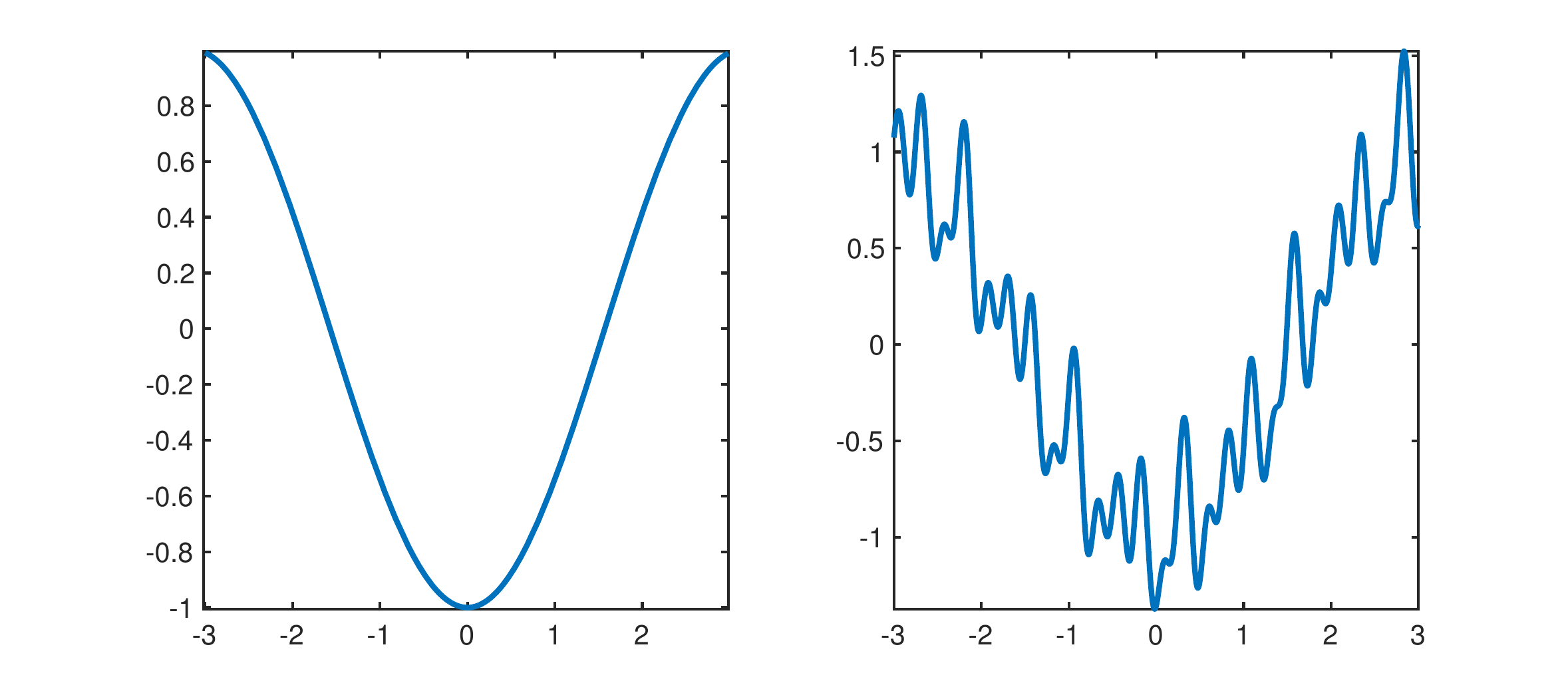}
	\caption{Two functions with the same global minimizer but different total variations. {\it Left}: the convex function  $f(x)=-\cos(x)$ defined over $(-3, 3)$; {\it Right}: the function added with a few high-frequency cosine modes.}
	\label{fig:roughness-cartoon}
\end{figure}

The  definition of TV is free of the deformation in the input variable: let $\varphi:[a',b']\to [a,b]$ be a diffeomorphism,
then $\TV(f\circ \varphi)=\TV(f)$.
For two functions defined on the same domain and have the similar size of the range, 
the TV norm can effectively describe the heuristic concept of ``roughness''. 
Refer to Figure \ref{fig:roughness-cartoon} where the right-side function has a much larger TV.
%fplot(@(x)  -cos(x),[-3,3]); axis([-3 3 -1 1])
% fplot(@(x)  -cos(x)-cos(x.*10)./5-sin(x.*15+0.4)./5+sin(x.*25-0.4)./5,[-3,3]);
%set(gca,'FontSize',20,'FontWeight','normal','LineWidth',2)
If $f$ is monotonic,  then $\operatorname{TV}(f)=\max  f  - \min f$.
There is another important  interpretation for the difference in the two functions in Figure \ref{fig:roughness-cartoon}
from the viewpoint of SGD
\cite{pmlr-v80-daneshmand18a,pmlr-v80-kleinberg18a}.
If one applies the SGD to minimize these two functions, it   takes 
much more time on the ``more rough'' function to reach the (global) optimal solution near $x=0$; the momentum acceleration like Adam
can mildly mitigate this slow convergence but generally speaking, the   function with a larger  TV is indeed harder to train.
Of course, the full gradient method without noise injection fails to obtain the global minimum for the non-convex function   in this case.
The above interpretation  of using the TV   to  describe the impact to  the stochastic training method can be explained
more precisely  from
the perspective  of the the well-known 
Freidlin-Wentzell large deviation theory \cite{FW2012,LTE2017SME,hu2019quasi}
for $$\mathrm{d}X_t = -\nabla f(X_t)+\sqrt{2\epsilon}\mathrm{d}W_t.$$
In this theory, the probability for the trajectories $X_t$ between two given endpoints are approximately (up to the exponential scale) determined by the 
so-called  {\it quasi-potential} function, 
for small $\epsilon$. We refer to the global minimum point in Figure \ref{fig:roughness-cartoon} as $o$.
Then  the quasi-potential $Q(o\to a)$ for transition starting from the lowest point $o$ and exiting the domain through the endpoint  $a$,
is the sum of all energy barriers \footnote{The barrier is the difference in $f$ between a local minimizer and its neighboring saddle point along the transition path.}. 
Therefore  we have  $\sum_{i=a,b}Q(o\to i)+ Q(i \to o)=\TV(f)$ holds {\it exactly} for any 1D  function defined over $[a,b]$. In this sense,
$\TV(f)$ represents how difficult   the stochastic gradient descent   approaches the lowest point $o$ from one boundary of the domain
{\it and} then exits the domain via either of  boundary points.
The bound of TV is also closely relevant to the magnitudes of the Fourier coefficients.
It is well known that a large Fourier coefficient at high frequency implies the function in space is more ``oscillatory''.
If $f$ on $[-\pi,\pi]$ has a bounded TV, then its Fourier coefficients $\hat{f}_{k}$ decay at least $O(1/k)$:  specifically we have \cite{bachmann2012fourier}:
$$\abs{\hat{f}_{k}}\leq \frac{2}{k\pi} \TV(f).$$
A small $\TV(f)$ corresponds to  small Fourier coefficients.

It is easy to see that $\TV(\alpha f) =\alpha  \TV(f), \, \alpha>0$. But 
to minimize $f$ and $\alpha f$ is exactly the same computational tasks if  the learning rate is rescaled accordingly. 
So, the index
for the function should be free of such dilation operation, 
and as a result  we propose the following modified TV
\begin{equation} \label{TT}
{T}(f):=\frac{1}{b-a}\frac{1}{[f]}\TV(f)=\frac{1}{b-a}\frac{1}{[f]}\int_a^b |f'(x)|\mathrm{d}x,
\end{equation}
where $$\displaystyle [f]=\max_{a\le x\le b} f(x)-\min_{a\le x\le b} f(x).$$
The denominators in \eqref{TT} for  the   domain size and range size
rescale  the graph of the function to ``fit'' into a unit square. 

Without loss of generality, we make the interval symmetric around the origin: $a=-b$.
Then if let $g(x)=\alpha f(\beta x)$ with two scalars $\alpha,\beta>0$ defined on the interval $[ a/\beta, b/\beta ]$, 
one can verify that $\TV(g)=\alpha \TV(f)$, but 
$T(g)=  \beta  T(f)$ due to the change of the interval size,  which suggests  an increasing roughness if $\beta$ is bigger than one
and this index $T$ is insensitive to $\alpha$. 
When $\beta>1$ and is an integer, by periodically extending the definition of $f$, 
we  now regard $g(x)=\alpha f(\beta x)$ defined on the same  $[a,b]$
as the original $f$ --- a conventional setting in homogenization theory\cite{papanicolau1978asymptotic}.
Then $\TV(g)=\alpha \beta \TV(f)$ and we still have $T(g)=\beta T(f)$ again since $\alpha$
is absorbed by the rescaling factor $[f]$ in the definition of \eqref{TT}.
One more property of $T$ is the following.
Assume  $f$  is  an even function attaining the minimum zero value  at  the origin in the interval  $I=[a,b]=[-l,l]$, then 
if $f$ is convex (or concave), we have $\TV(f)=2[f]$, and $T(f)\equiv 1/l$.  One example like this  is 
the quadratic function $f(x)=\beta x^2/2$. If $f$ is not even, then  $T$ in \eqref{TT} is sensitive to
the values at two endpoints.

\subsection{Roughness index for high dimensional functions}
To generalize the above 1D index $T$ to any dimension, 
we follow the idea of projection to randomly sampled direction with filter-wise normalization
$$f_d(s):= \mathcal{J}(\theta+ sd)$$
where $\theta$ is a given reference point and $d$ is a Gaussian random direction with zero mean and identity covariance matrix followed by  filter-wise normalization \cite{li2018visualizing}. 
The domain of $s$ is defined on  a prescribed interval $[-l,l]$. By varying $l$, we can change the size of the region in concern around the reference point $\theta$.
Unlike in \cite{li2018visualizing} 
which used just one sampled direction $d$ in the visualization procedure, 
we consider the standard deviation  of $f_d$ with respect to the randomness in the directions, so the {\bf roughness index} ({\bf RI}) is defined as follows
\begin{equation}\label{index}
\mathcal{I}(\mathcal{J};\theta) = \frac{\mathbf{std}_d  T(f_d)} {\mathbf{E}_{d}T(f_d)} .
\end{equation}
Here the standard deviation is adopted to describe the 
change of ``roughness''  across different directions. 
The rescaling  by the  expectation here is to further reduce the influence of the magnitude of $T$ values.

\begin{example} \label{ex:quad}
	We examine the index  by looking at  a quadratic landscape  $\mathcal{J}(\theta) = \frac12 \theta^\top H \theta$ where the reference point is taken as the minimizer (the origin)
	and set the interval size $l=1$.
	$H$ is a positive definite matrix.
	Then  $f'_d(s)=s d^\top H d$ and $\TV(f_d) = \abs{d^\top H d}$. 
	If $d$ follows the standard Gaussian distribution with zero mean and identity covariance matrix, 
	then   by Hutchinson's trick, $\mathbf{E}_d \TV(f_d)= \mathbf{E}_d d^\top H d =\mathbf{E}_d  \Tr(dd^\top H)  ={\Tr}( \mathbf{E}_d  (dd^\top)  H)=\Tr(H)$.
	But $T(f_d)\equiv 1$ in view of \eqref{TT} and the roughness index $\mathcal{I}$ in \eqref{index} is zero for   {any} quadratic function.
\end{example}

\subsection{Algorithm}
The details of the computational procedure is as follows. Assume $\theta^*$ is an arbitrary  point of interest.
In many cases, we consider a minimum point obtained by minimizing  the loss function $\mathcal{J}$. To calculate RI  w.r.t. this point, detailed description on the numerical implementation of RI is available in Algorithm 1.
The complexity is linearly proportional to $M\times m$ and independent of the dimension of $\theta$.

\begin{algorithm}
	\caption{{\bf Computation of Roughness Index}}
	\KwIn{Loss $\mathcal{J}$, point  $\theta^*$, number of directions $M$,   interval length $l_i$  and  number of step size $m$}
	\KwOut{Roughness Index $\mathcal{I}$ at  $\theta^{*}$}

% 	\KwIn{Loss $\mathcal{J}$, point of interest $\theta^*$, number of iid Gaussian random direction $M$, size of interval $l_i$(in practice, we use the same $l$ for all $l_i$ in all directions) and grid size of interval $m$}
% 	\KwOut{Roughness Index $\mathcal{I}$ at  $\theta^{*}$}
	
	$i \longleftarrow  1$\\
	
	\While{$i \le M$}{
		Sample an iid standard Gaussian random direction $d_i$; \\
		Apply the filter-wise normalization for $d_i$: $\bar{d_i} \gets d_i$ \\
		$j \longleftarrow  0$\\
		\While{$j \le m$}{
			Partition   $[-l_i, l_i]$ into $m + 1$ subintervals uniformly:
			$$s_{i, j} = -l_i + j \frac{2l_i}{m}, j = 0, 1, \cdots, m$$ \\
			$j \longleftarrow  j + 1$
		}
		Calculate the maximum and minimum along  $\bar{d_i}$: \\
		$$\mathcal{J}_{\max}^i = \max_{0 \le j \le m}  \{\mathcal{J}(\theta^* + s_{i, j} \bar{d_i})\}$$
		$$\mathcal{J}_{\min}^i = \min_{0 \le j \le m} \{\mathcal{J}(\theta^* + s_{i, j} \bar{d_i})\}$$
		
		Approximate  normalized TV $T_i$ : \\	
		
		$$T_i =  \frac{1}{2l_i} \sum_{j = 0}^{m-1}  \frac{|{\mathcal{J}(\theta^* + s_{i, j} \bar{d_i}) - \mathcal{J}(\theta^* + s_{i, j+1} \bar{d_i})}|} {\mathcal{J}_{\max}^i - \mathcal{J}_{\min}^i}$$
		
		$i \longleftarrow  i + 1$\\
	}
	
	The roughness index $\mathcal{I} := {\sigma}/{\mu}$, where $\mu,\sigma$ are the mean value and the standard deviation of ${\{T_i\}}_{i=1}^M$.

\end{algorithm}

The number of directions $M$ and the  number of partitions for interval $m$ are chosen sufficiently large in practice to make sure the numerical results are convergent. 
In  addition, the various values of  interval length $l$ are also tested for specific applications (See Remark \ref{rem1}).

\section{Numerical Results}\label{sec:result}

Consider the Poisson equation on $\Omega = {(0, 1)}^d$:
\begin{equation}\label{eqn: poisson equation}
\begin{cases}
- \Delta u = f(x), & \; \text{in} \; \Omega,\\
u(x) = 0, & \; \text{on} \; \partial \Omega.
\end{cases}
\end{equation}
The forcing term $f$ is specified by assuming the form of the   solution first.
For example, we assume    the exact solution 
\begin{equation} \label{uex}
u(x) = \prod_{i = 1}^{d} \sin (\pi x_i), \; x = (x_1, \cdots, x_d),
\end{equation}
then we  have $f(x) = d \pi^2 \prod_{i = 1}^{d} \sin (\pi x_i)$.
Denote 
\begin{equation}\label{neural solution}
u(x;\theta) =   \prod_{i = 1}^{d} (x_i - 1) x_i \cdot \NN(x;\theta).
\end{equation}
where $\NN(x;\theta)$ is a function represented by a NN. The corresponding loss functions are
\begin{equation}\label{eqn:lossdgm_d}
\mathcal{J}_{\textrm{G}} (\theta) =  \int_{\Omega}\left(- \Delta u(x;\theta) - f(x)\right)^2\mathrm{d} x
\end{equation}
for the DGM, and
\begin{equation}\label{eqn:lossdrm_d}
\mathcal{J}_{\textrm{R}} (\theta) =  \int_{\Omega} \left(\frac{1}{2} {|\nabla u(x;\theta)|}^2 - f(x) u(x;\theta) \right)\mathrm{d} x
\end{equation}
for the DRM, respectively. 
% $\nabla$  refers to   the gradient  w.r.t. $x$ and $|\cdot|$ is the Euclidean norm in $\mathbb{R}^{d}$.
% Note that in this case,  the loss functions $\mathcal{J}(\theta)$ are generally non-convex.

In what follows, we use the relative $L^2$ error to measure the  numerical error
of solving the PDE,
\begin{equation} \label{def:err}
\textrm{error} = \frac{{\begin{Vmatrix}u(x;\theta^*) - u(x)\end{Vmatrix}}}{{\begin{Vmatrix}u(x)\end{Vmatrix}}},
\end{equation}
where $\begin{Vmatrix}\cdot\end{Vmatrix}$ denotes the $L^2$ norm for functions of $x$, $u(x;\theta^*)$ is the DNN approximation, and $u(x)$ is the exact solution. 

\subsection{1D Poisson equation}  Consider the following 1D Poisson equation
\begin{equation}\label{eqn: poisson equation 1d}
\begin{cases}
- u''(x) = f(x), & \; x\in (0, 1),\\
u(0) = u(1) = 0. & 
\end{cases}
\end{equation}
The exact solution is set as $u(x)=\sin \pi x $, so that $f(x)=\pi^{2} \sin \pi x$.
At this  true solution,
 we have  the global minima for  $\mathcal{J}_{\textrm{G}} (u(x)) = 0$, and $\mathcal{J}_{\textrm{R}} (u(x)) = -\pi^2 / 4 \approx -2.4674$.

The numerical solution is in the form of
$
u(x;\theta) =    (x  - 1) x  \cdot \NN(x;\theta).
$
% We  use  the   neural network as specified in Section \ref{ssec:nn}.
Various width $w$ is tested for ResNet and 
FCNet.
The loss  functions $\mathcal{J}(\theta)$ are non-convex now,
but in practice one can generally find the global minima
due to the perfect fitting capability of the neural network \cite{gamst2017energy,45820}.

The 1D integrals in \eqref{eqn:lossdgm_d} and \eqref{eqn:lossdrm_d}  are approximated by a 
quadrature rule with  $N$ uniform points on the interval $[0,1]$. And we refer this $N$ as to the batch size since in the training
we use all these $N$ points in each gradient-based iteration.

\subsubsection{ Local minimizers}

\begin{table}[htbp!]
	\centering
	\caption{Losses at  $\theta_G, \theta_R$ and $ \wt{\theta}_G$.
	The   (global) minimum values of $\JG$ and $\JR$ are $0$ and $-\frac{\pi^2}{4}\approx -2.4674$.
		We treat $\theta_G$ and $\wt{\theta}_G$ as the two local 
		minimizers of $\JG$ and all three    as local 
		minimizers of $\JR$.}
	\label{tab: threeloss}
	\begin{tabular}{c|c|c|c}
		\toprule[1pt]
		
		\noalign{\smallskip}
		\multirow{1}*{loss} 
		&\multicolumn{1}{c|}{$\theta_G$} & {$\theta_R$} & {$ \wt{\theta}_G$}  \\		
		\noalign{\smallskip}
		\midrule[1pt]
		\noalign{\smallskip}
		\multirow{1}*{$\mathcal{J}_G(\theta)$}
		& $5.9933\text{e-05}$ & $0.1044$ & $5.7418\text{e-05}$  \\
		\noalign{\smallskip}
		
		\multirow{1}*{$\mathcal{J}_R(\theta)$}
		& $-2.4715$ & $-2.4716$ & $-2.4715$ \\
		\noalign{\smallskip}
		\bottomrule[1pt]
	\end{tabular}

\end{table}

\begin{table}[htbp!]
	\centering
	\caption{The distance between  $\theta_G, \theta_R$ and $ \wt{\theta}_G$.}
	\label{tab: threedistance}
	\begin{tabular}{c|c|c|c}
		\toprule[1pt]
		
		\noalign{\smallskip}
		\multirow{1}*{distance} 
		&\multicolumn{1}{c|}{$(\theta_G, \theta_R)$} & {$(\theta_G, \wt{ \theta}_G)$} & {$(\theta_R, \wt{ \theta}_G)$}  \\		
		\noalign{\smallskip}
		\midrule[1pt]
		\noalign{\smallskip}
		\multirow{1}*{$\begin{Vmatrix} \cdot \end{Vmatrix}_2 $}
		& $3.7243$ & $3.8342$ & $0.3349$  \\
		\noalign{\smallskip}
		
		\multirow{1}*{$\begin{Vmatrix} \cdot \end{Vmatrix}_\infty$}
		& $2.2392$ & $2.3052$ & $0.2138$ \\
		\noalign{\smallskip}
		\bottomrule[1pt]
	\end{tabular}

\end{table}

Staring  from the {\it same} initial guesses used to train $\JG$ and $\JR$,
we use  the full-batch gradient descent to 
 find one local minimizer for each loss function, denoted by 
$\theta_G$ and $\theta_R$, respectively.
Even though both parameters $\theta_G$ and $\theta_R$ gives approximate solutions to the PDE,
these two  parameters $\theta_G$ and $\theta_R$ are quite different. See Table \ref{tab: threedistance}.
After  obtaining  $ {\theta}_G$ and $ {\theta}_R$   from the DGM and DRM respectively, 
we swap them as the new initial guesses to train $\JG$ and $\JR$.
This is to look for  a new optimal  parameter $\wt{\theta}_G$    by minimizing $\JG$ with the 
	new initial guess $\theta_R$ and   for $\wt{\theta}_R$  of  $\JR$ in a like manner by using the initial  $\theta_G$.
We find that   $\wt{\theta}_R$ is almost identical to $ {\theta}_G$ and conclude $\theta_G$ and $\wt{\theta}_G$ are   minimizers of $\JG$; 
	$\theta_R$ and ${\theta}_G$ (  $=\wt{\theta}_R$)  
	as well as $\wt{\theta}_G$ are   minimizers of $\JR$.
	The  loss values at these points are shown in Table \ref{tab: threeloss}.

  \subsubsection{Difference between DGM and DRM}
We observed that the DGM generally obtains a better  accuracy in {\it solving PDE}  result than the DRM in our case here.
We  compare  their accuracy by
checking the PDE errors in \eqref{def:err} of their corresponding PDE solutions  $u(\cdot;\theta_G)$
and $u(\cdot;\theta_R)$.
	We tested the ResNet of one block with different widths
	in Table \ref{tab: Relative error for DGM and DRM}.
	  Since the   NN and the training algorithm as well as the initial guess are exactly the same, we attribute this discrepancy to the difference of loss in the DGM and DRM.
	\begin{table}[htbp!]
		\caption{  {The relative $L^2$ PDE error 
				defined by   \eqref{def:err} for deep Galerkin method and deep Ritz method after   training  10000 epochs with  different widths of the ResNet.}
		} \label{tab: Relative error for DGM and DRM}
		\centering
		\begin{tabular}{c|c|c|c|c|c}
			\toprule[1pt]
			
			\noalign{\smallskip}
			\multirow{1}*{$w$} 
			&\multicolumn{1}{c|}{$2$} & {$3$} & {$4$}  & {$5$} & {$6$}  \\		
			\noalign{\smallskip}
			\midrule[1pt]
			\noalign{\smallskip}
			\multirow{1}*{$u(\cdot;\theta_G)$}
			& $\text{5.21e-2}$ & $\text{1.81e-2}$ & $\text{7.12e-4}$  & $\text{8.01e-8}$ & $\text{8.31e-8}$\\
			\noalign{\smallskip}
			\multirow{1}*{$u(\cdot;\theta_R)$}
			& $\text{1.64e-3}$ & $\text{9.48e-4}$ & $\text{7.63e-4}$  & $\text{7.76e-4}$ & $\text{6.75e-6}$\\
			\noalign{\smallskip}
			\bottomrule[1pt]
		\end{tabular}

	\end{table}

	\begin{figure}[htpb!]
		\centering

		\subfigure[The decay of loss functions. ]{
			\includegraphics[width=1.8in]{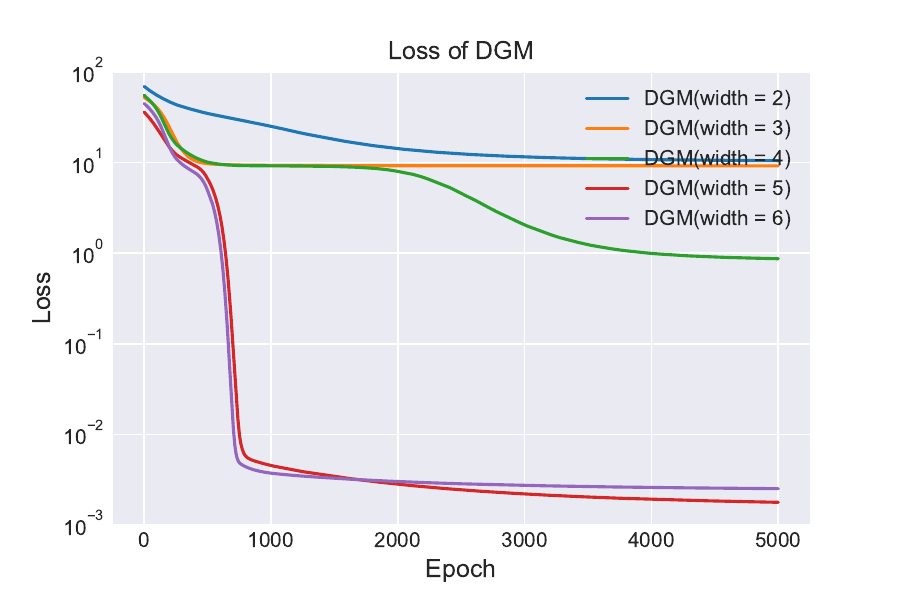}
			\includegraphics[width=1.8in]{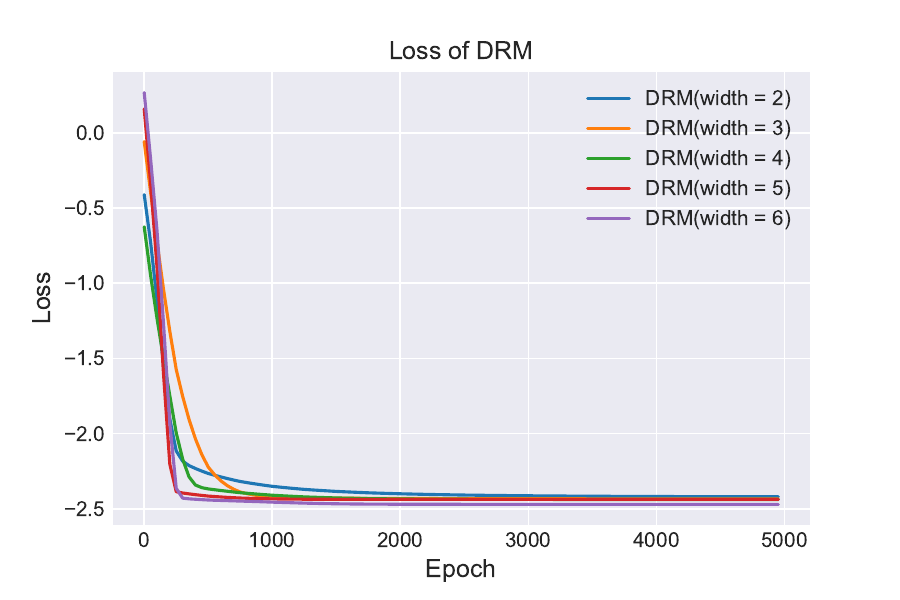}
		}
		\subfigure[The decay  of relative $L_2$ error of $u(x;\theta)$ to the true PDE solution.  ]{
			\includegraphics[width=1.8in]{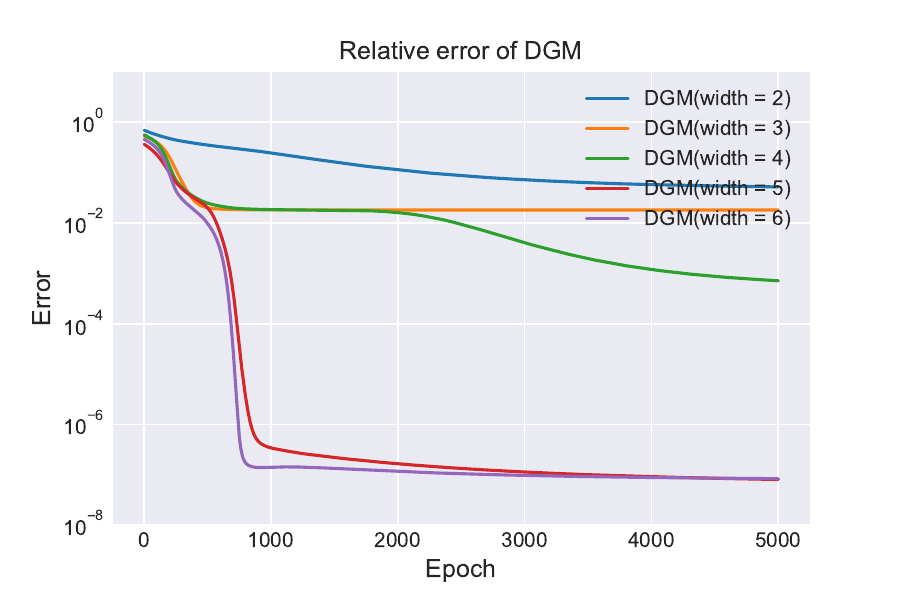}
			\includegraphics[width=1.8in]{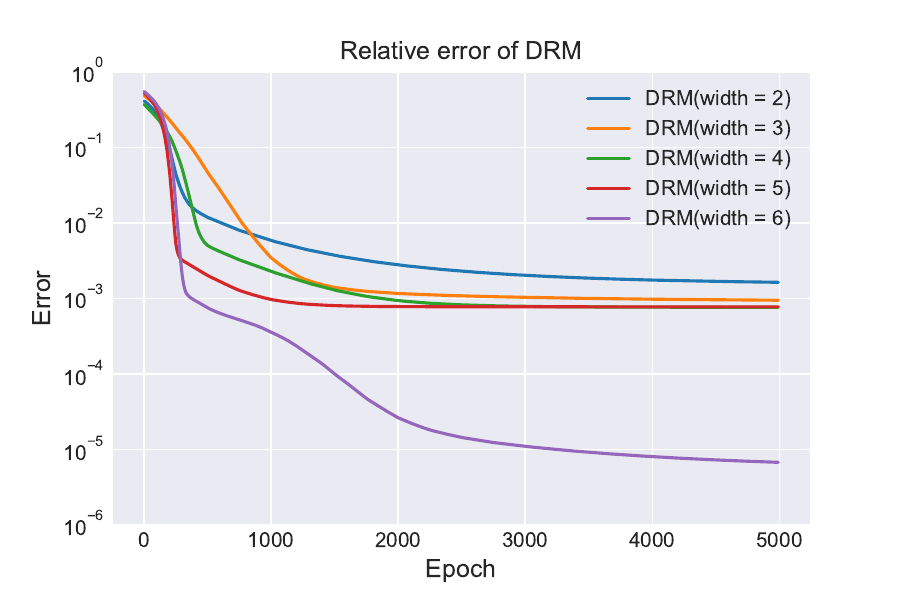}
		}
				\caption{The loss functions and  the relative $L^2$ error for   ResNet  with   width $w = 2, 3, 4, 5, 6$.  {\it Left column}: DGM; {\it Right column}: DRM. }
		\label{fig: Loss and error for widenn}

	\end{figure}
	We furthermore  provide  complementary results about the convergence for DGM and DRM toward $\theta_{G}$ and $\theta_{R}$ respectively.	
	Figure \ref{fig: Loss and error for widenn}  shows the decay of the loss and the  relative $L^2$ error \eqref{def:err} in the training process.
		One interesting observation comes from the comparison of the loss and the  error.
	The DRM is very effective to decrease the loss for all widths, but inefficient to decrease the PDE error.
	It seems that after the early stage of quick decay for the  loss function, 
	the    DRM trajectories
  wander around in a neighbor of the   minimizer  of  the loss function  
	in order to further reduce the PDE's  error, but with much more strenuous efforts than the DGM.
	As a comparison, 
	the DGM has a better  match for 
	the decay between the PDE error  and   the loss function.
This is easy to understand since  	by \eqref{eqn:lossdgm_d}, 
	the loss of the DGM is $\mathcal{J}_{\textrm{G}} (u) = \int_0^1  ( u'' - u_{\textrm{ex}} '' )^2 \mathrm{d} x =\|u'' - u_{\textrm{ex}} ''\|^2,$
	with the only difference of a (linear) Laplace operator, 
	which is  more  closely linked  to the PDE error \eqref{def:err} than the DRM.

	\subsubsection{Roughness index (RI)} 
	Now we report our  main  numerical results of $\mathcal{I}$ for this 1D problem.  
	We  record roughness indices  in several setting of parameter combinations. The calculation  involves the minimizers  of interests, the number of  directions $M$, the   interval length  $l$,
	the number of points  $m$ partitioned in the interval. 
	
	We first present the results of roughness indices of the DGM and the DRM 
	around their first set of optimal parameters $\theta_G$ and $\theta_R$.
	With a fixed width $w=4$, Table \ref{tab: Total variation compara_M} to 
	Table \ref{tab: tvd for widenn} show the comparing results of the roughness indices for the two models with  various combinations of network architecture (ResNet or FCNet), the width $w$, the values of $M$, $l$ and $m$.
	In all cases, particularly with the ResNet architecture, we have strong numerical evidences to claim that {\it the roughness index in the DGM is significantly smaller than that in the DRM}.

	\begin{table}[tbhp!]
		\centering
			\caption{RI for different $M$ with $l= 0.0001$ and $m = 30$.}\label{tab: Total variation compara_M}
		\begin{tabular}{c|cc|cc}
			\toprule[1pt]
			
			\noalign{\smallskip}
			\multirow{2}*{$M$} 
			&\multicolumn{2}{c|}{$\mathcal{I} _{DGM}$} &\multicolumn{2}{c}{$\mathcal{I} _{DRM}$}\\			
			&\multicolumn{1}{c}{ResNet} & \multicolumn{1}{c|}{FCNet} & \multicolumn{1}{c}{ResNet} & \multicolumn{1}{c}{FCNet}\\
			\noalign{\smallskip}
			\midrule[1pt]
			\noalign{\smallskip}
			\multirow{1}*{50}
			& $0.0455$ & $0.2387$  & $0.4665$ &$ 0.2472$\\
			\multirow{1}*{100} 
			& $0.0615$ &$ 0.2157$ & $0.4443$ & $0.2256$\\
			\multirow{1}*{150} 
			& $0.0668$ &$ 0.2186$ & $0.4653$ & $0.2195$\\
			\noalign{\smallskip}
			\bottomrule[1pt]
		\end{tabular}

	\end{table}

	\begin{table}[htbp!]
		\centering
		\caption{RI for  different $l$ with $M = 100$ and $m = 100$.}\label{tab: Total variation compara_l_i}
		\begin{tabular}{c|cc|cc}
			\toprule[1pt]
			
			\noalign{\smallskip}
			\multirow{2}*{$l$} 
			&\multicolumn{2}{c|}{$\mathcal{I} _{DGM}$} &\multicolumn{2}{c}{$\mathcal{I} _{DRM}$}\\			
			&\multicolumn{1}{l}{ResNet} & \multicolumn{1}{c|}{FCNet} & \multicolumn{1}{c}{ResNet} & \multicolumn{1}{c}{FCNet}\\
			\noalign{\smallskip}
			\midrule[1pt]
			\noalign{\smallskip}
			\multirow{1}*{0.00025}
			& $0.0287$ & $0.1846$ & $0.6743$ & $0.2139$\\
			\multirow{1}*{0.0005} 
			& $0.0073$ & $0.1336$ & $0.7264$ & $0.1712$\\
			\multirow{1}*{0.001} 
			& $0.0109$ & $0.0731$ & $0.7311$ & $0.1291$\\
			\multirow{1}*{0.005}
			& $0.0074$ & $0.0253$ & $0.1863$ & $0.0537$\\
			\multirow{1}*{0.01} 
			& $0.0127$ & $0.0157$ & $0.1525$ & $0.0227$\\
			\multirow{1}*{0.05} 
			& $0.0418$ & $0.0553$ & $0.0876$ & $0.0705$\\
			\noalign{\smallskip}
			\bottomrule[1pt]
		\end{tabular}

	\end{table}

	\begin{table}[htpb!]
		\centering
		\caption{RI  for different  $l$ and $m$.  ($M = 100$ and   ResNet.)}\label{tab: Total variation compara_l_i_m}
		
		\begin{tabular}{c|c|c|c}
			\toprule[1pt]
			$l$ & $m$  & $\mathcal{I}_{DGM}$ & $\mathcal{I}_{DRM}$ \\
			\midrule[1pt]
			0.00005 &  20  & $0.0517$ & $0.3639$ \\
			
			%		0.00010 &  30  & $  {0.0615}$ & $ {0.4443}$ \\
			%		\hline
			%		0.00010 &  40  &  $0.0576$ & $0.4542$ \\
			%		\hline
			0.00010 &  50  & $0.0587$ & $0.4593$ \\
			
			0.00015 &  60  & $0.0394$ & $0.5709$ \\
			
			0.00020 &  80  & $0.0353$ & $0.6222$ \\
			
			0.00025 &  100  & $0.0287$ & $0.6743$ \\
			
			0.00030 &  120  & $0.0275$ & $0.7096$ \\
			\bottomrule[1pt]
		\end{tabular}

	\end{table}

	\begin{table}[htpb!]
		\centering
			\caption{RI for neural networks with  width $w = 2, 3, 4, 5, 6$. ($l = 0.02$, $M = 100$, $m = 100$, and the ResNet.)}\label{tab: tvd for widenn}
		\begin{tabular}{c|c|c}
			\toprule[1pt]
			$w$ & $\mathcal{I}_{DGM}$ & $\mathcal{I}_{DRM}$ \\
			\midrule[1pt]
			2 & $0.0356$ & $0.0843$ \\
			
			3 & $0.0289$ & $0.2389$ \\
			
			4 & $0.0216$ & $0.0890$ \\
			
			5 & $0.0266$ & $0.0992$ \\
			
			6 & $0.0208$ & $0.0481$ \\
			\bottomrule[1pt]
		\end{tabular}
	
	\end{table}

	\begin{remark}\label{rem1}
		We remark that although the choice  of $M$ and $m$ is simple (the larger the better), the choice of the interval length $l$ is  important and one should test a few values for this parameter.
		$l$ characterizes the size of a small neighborhood we are interested when measuring  the roughness.   If $l$ is too large,  the domain of interest is too large to smear the 
		roughness around the reference point.     Table \ref{tab: Total variation compara_l_i_m_lager}  shows such  phenomena  as $l$ increases to a very large value:
		the disparity in the roughness index  between the two models is less and less significant. 
		The visualization plot in Figure \ref{fig: Loss landscape for FCNet and ResNet} corresponds to $l=0.01$.
		Conceptually, the suitable size of $l$ should be comparable to the size of the basin of attraction, 
		but here we deal with a highly non-convex landscape and it is not possible to pinpoint  this value. 
		So instead, we varied the choices of $l$ in practice and seek for a robust result in a reasonable range of $l$. We find $l=0.01$ is quite representative for our example here.
		
	\end{remark} 
	\begin{table}[htp]
		\centering
		\caption{Roughness index  tends to the same  for very large values of $l$. $M = 100$, $m = 100$, and the ResNet with $w=4$.}\label{tab: Total variation compara_l_i_m_lager}
		\begin{tabular}{c|c|c}
			\toprule[1pt]
			$l$ & $\mathcal{I}_{DGM}$ & $\mathcal{I}_{DRM}$ \\
			\midrule[1pt]
			0.1 & $0.0759$ & $0.1050$ \\
			
			0.2 & $0.1151$ & $0.1381$ \\
			
			0.3 & $0.1509$ & $0.1727$ \\
			
			0.4 & $0.1562$ & $0.1823$ \\
			\bottomrule[1pt]
		\end{tabular}
		
	\end{table}
	
	Lastly, we report the RI for the second set of parameters $\wt{\theta}_G$.
	Recall that 
	we validated  $\theta_G$ and $\wt{\theta}_G$ are two different minimizers of $\JG$; 
	$\theta_R$ and ${\theta}_G=\wt{\theta}_R$   are two different minimizers of $\JR$.
	We have reported the roughness index for $\theta_G$ and $\theta_R$ before.
	Table  \ref{tab: Index compara 1d}  adds the RI
	of the DGM and DRM at all these three points.
	It shows that the RI of the DGM is almost
	equal for the DGM's two local minimizers and this is also true for the 
	DRM's   local minimizers. And the roughness index of the DGM 
	is indeed much smaller than the roughness index of the DRM,
	regardless of which minimizer  of their own is investigated. 
	We can not confirm that this holds for all local minimizers 
	since it is not possible to explore all these minimizers. 
	But we are inclined to  the conjecture of a larger 
	roughness index for the landscape of the DRM than  the DGM, when the ResNet  is used.
	\begin{table}[htbp!]
		\centering
		\caption{RI at  different    reference points with $M = 100, l = 0.01$, $m = 100$, the ResNet, and width $w=4$. ``$*$'': Note that $\theta_R$ is not 
			a numerical minimizer of $\JG$, which eventually evolves to $\wt{\theta}_G$ by gradient descent.}\label{tab: Index compara 1d}
		\begin{tabular}{c|c|c}
			\toprule[1pt]
			
			\noalign{\smallskip}
			\multirow{1}*{reference point} 
			&\multicolumn{1}{c|}{$\mathcal{I} _{DGM}$} &{$\mathcal{I} _{DRM}$}\\		
			\noalign{\smallskip}
			\midrule[1pt]
			\noalign{\smallskip}
			\multirow{1}*{$\theta_G$}
			& $0.0127$ & $0.1525$ \\
			\multirow{1}*{$\theta_R$} 
			& $0.1448*$ & $0.1732$ \\
			\multirow{1}*{$\wt{\theta}_G$} 
			& $0.0153$ & $ {0.1660}$ \\
			\noalign{\smallskip}
			\bottomrule[1pt]
		\end{tabular}

	\end{table}

	\subsubsection{Validation of  RI  by visualization} 
	After  we calculated   the numerical values of RI for the DGM and  DRM models,
	we have reached a    conclusion that the landscape of the DRM 
	seems rougher than the DGM.
	To validate this claim, we apply the  visualization technique in  \cite{li2018visualizing}
   to show   heuristic and visual  evidence.

\begin{figure}[htbp!]
		\centering
		\subfigure[FCNet.   {\it Left  }: $\JG(\cdot)$ near 
		$\theta_G$; {\it Right  }:  $\mathcal{J}_{R}(\cdot)$ near 
		$\theta_R$.]{
			\includegraphics[width=1.8in]{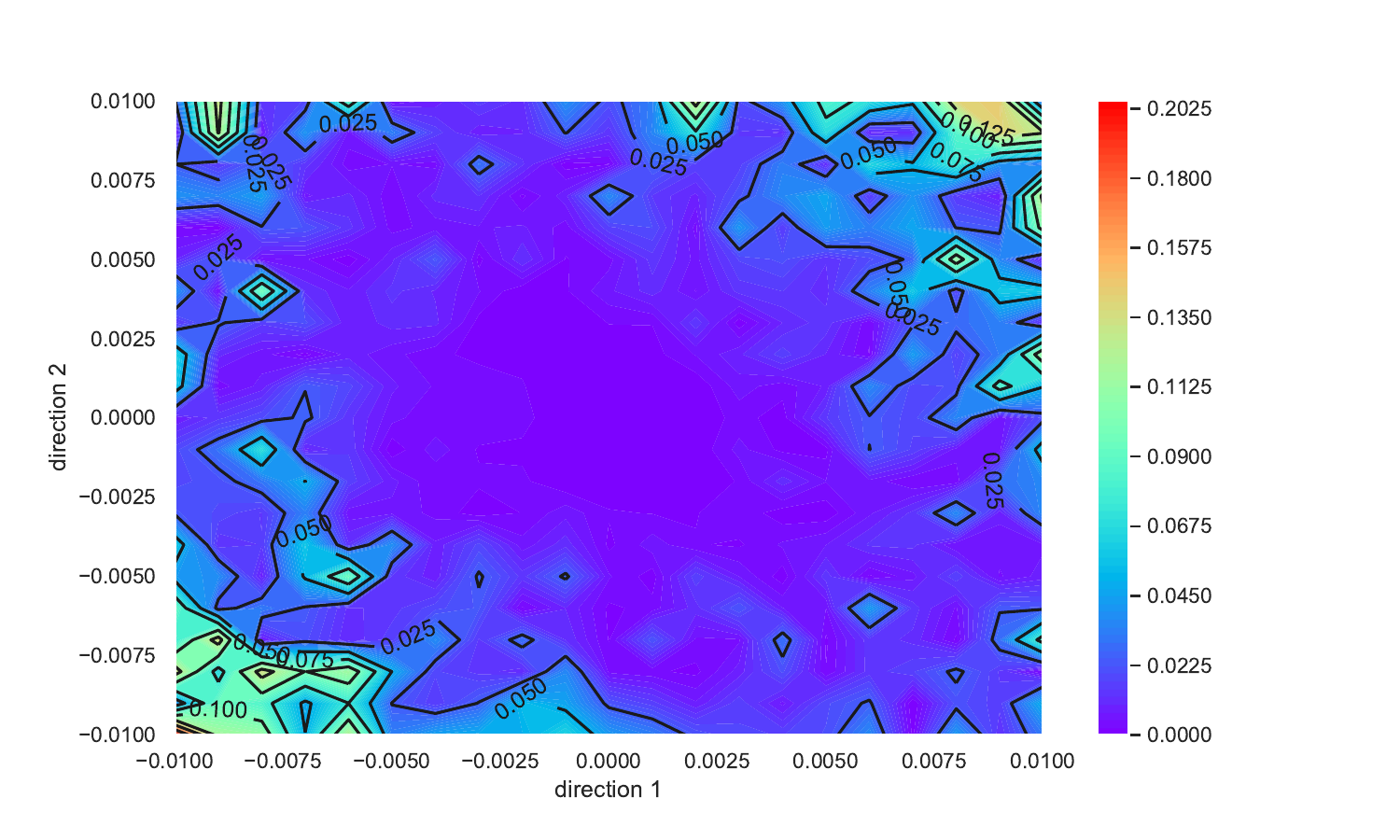}
			\includegraphics[width=1.8in]{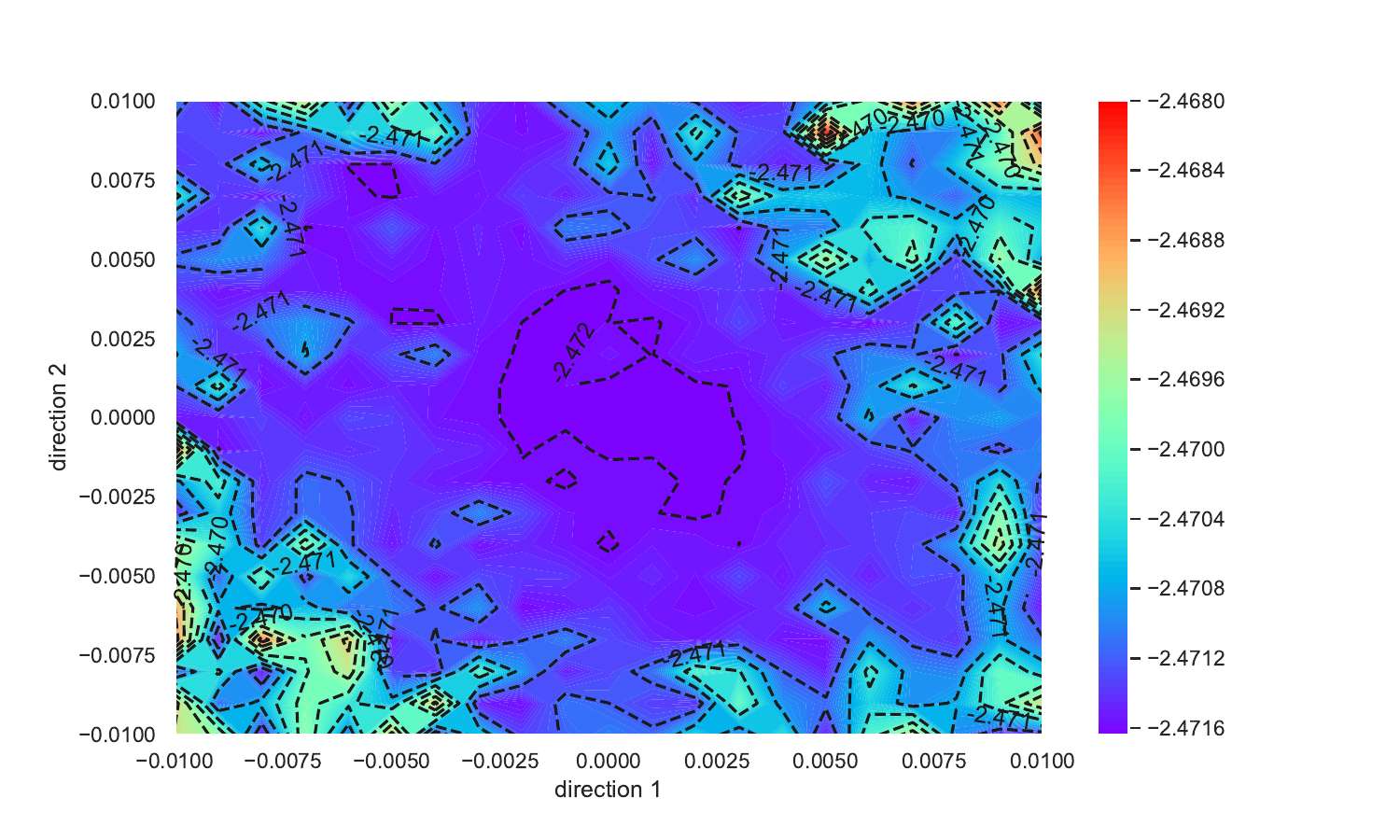}
		}  
		\subfigure[ResNet.   {\it Left  }: $\JG(\cdot)$ near   
		$\theta_G$; {\it Right  }:  $\mathcal{J}_{R}(\cdot)$ near 
		$\theta_R$.]{
			\includegraphics[width=1.8in]{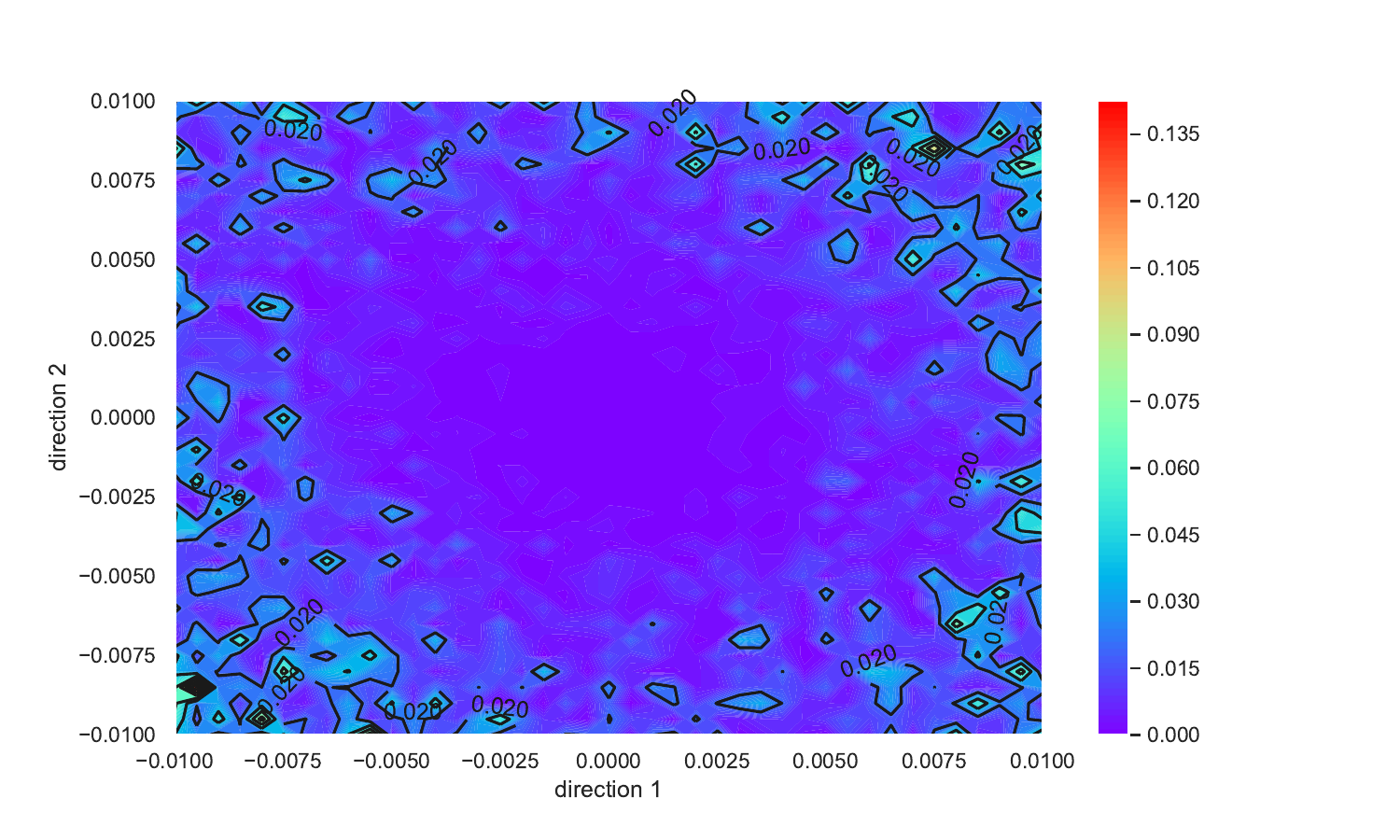}
			\includegraphics[width=1.8in]{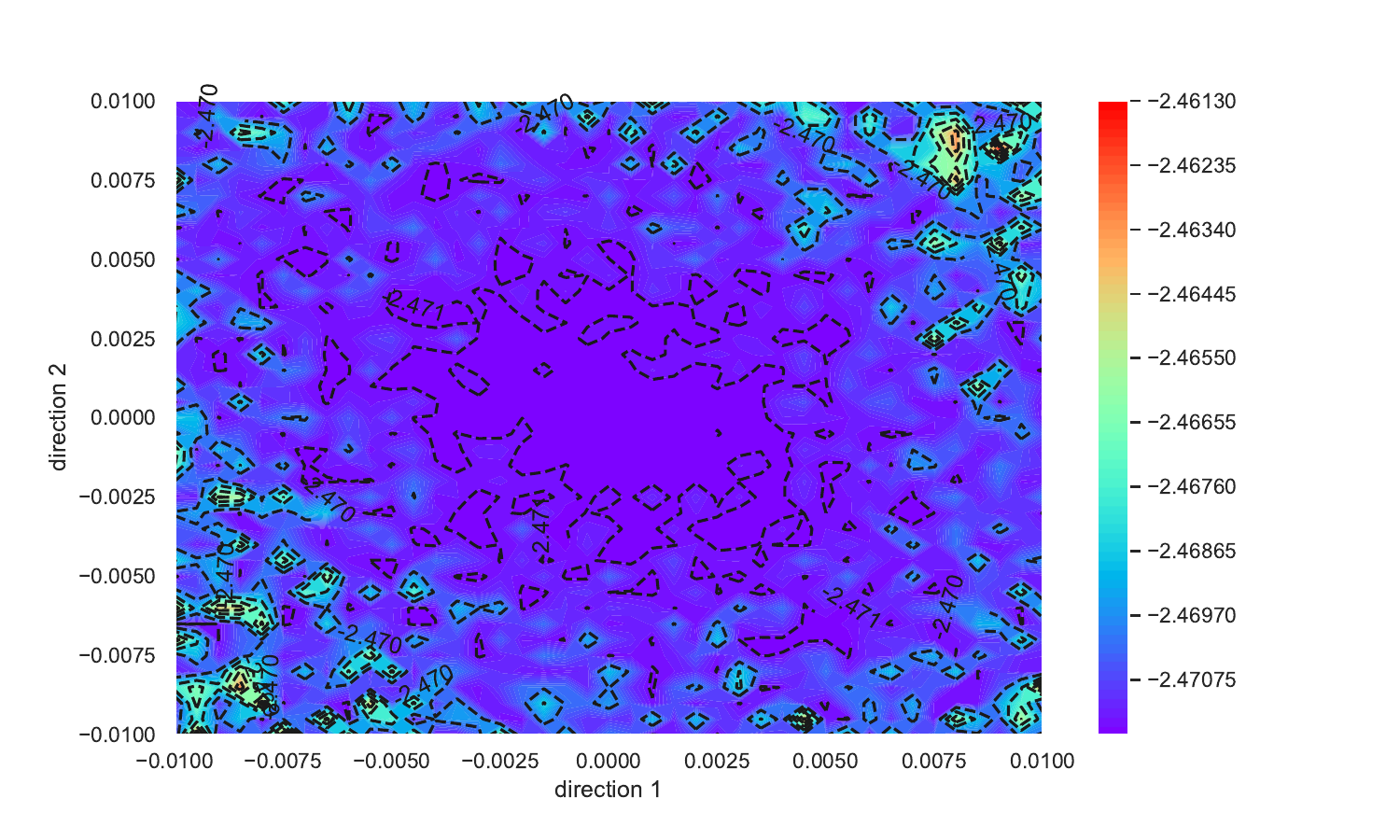}
		}  
		\subfigure[ResNet.   {\it Left  }: $\JG(\cdot)$ near   
		$\wt{\theta}_G$; {\it Right  }:  $\mathcal{J}_{R}(\cdot)$ near 
		$\wt{\theta}_R(=\theta_G)$.]{
			\includegraphics[width=1.8in]{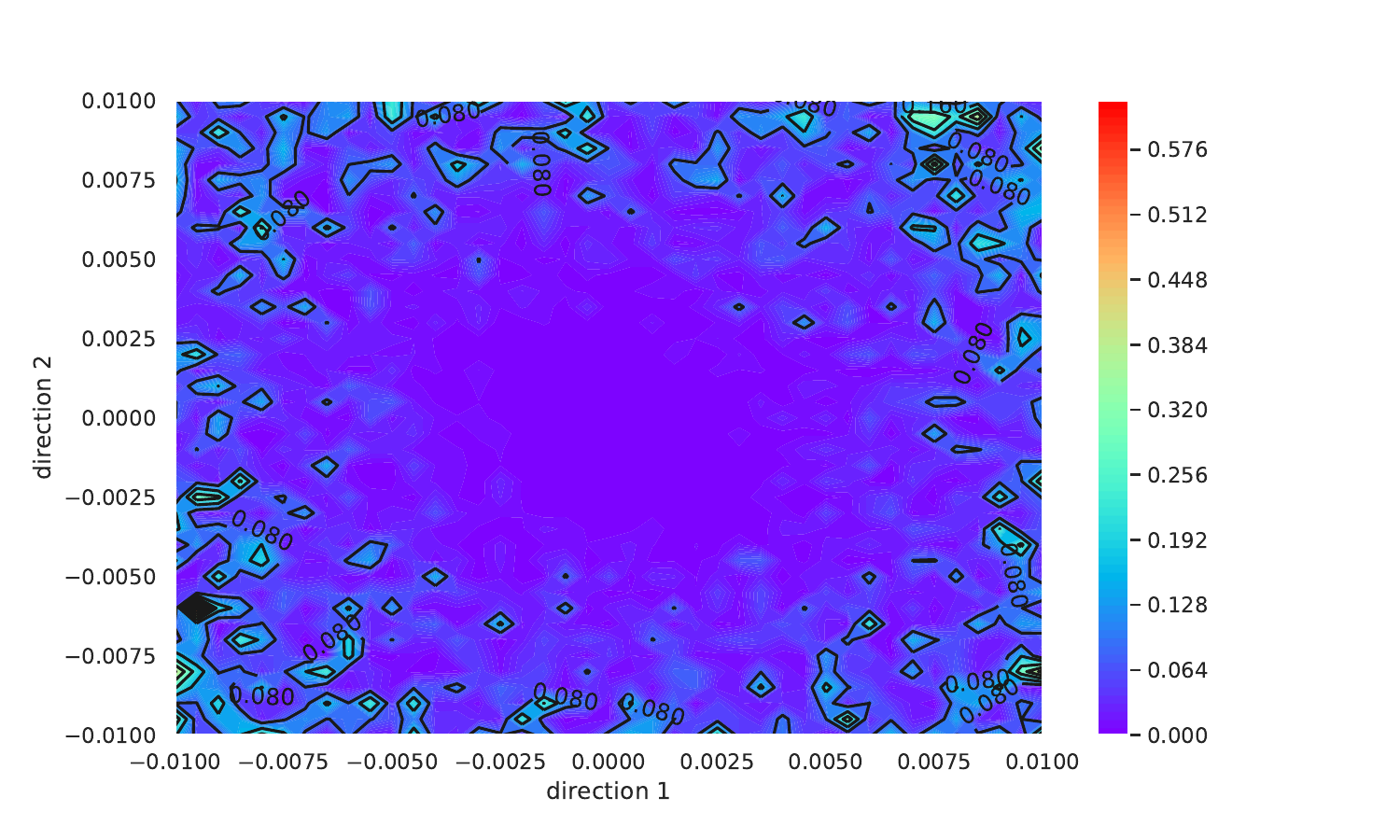}
			\includegraphics[width=1.8in]{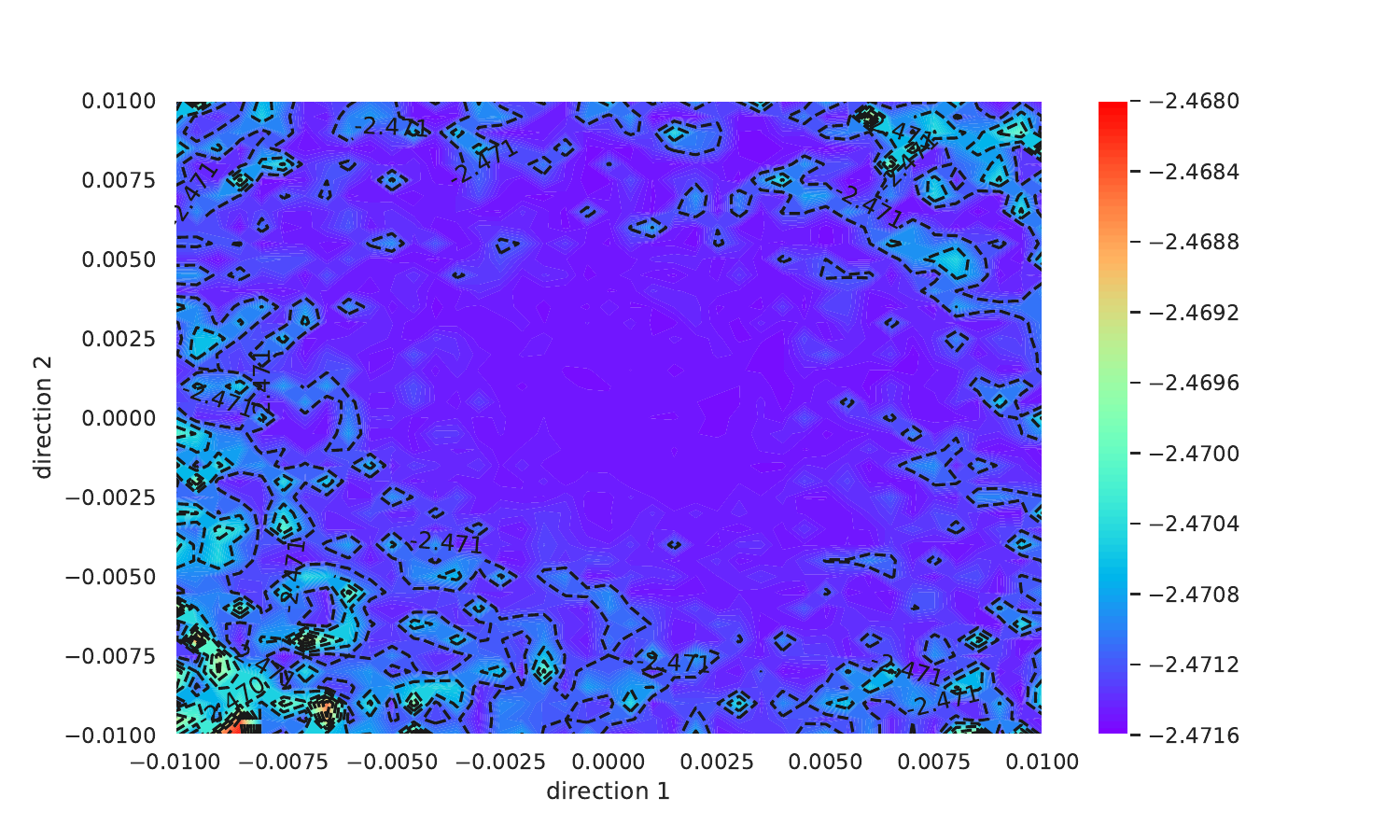}
		} 
		\caption{Contour plots for the two dimensional visualization  of loss landscapes around their   local minimizers.
			In each figure, the contour plot contains exactly  eight isolines   with  equal gaps between the minimal and the maximal values
			(marked    in the vertical colorbars).
			The left panels refer to the loss landscape
			$\JG$ while the right panels refer to the loss landscape	$\JR$. 
			(Batch size 200 and neural width $w = 4$.)
		}	\label{fig: Loss landscape for FCNet and ResNet}
	\end{figure}

	We use visualization with filter-wise normalization in a randomly chosen 2D space.
	The contour plots of loss landscapes  for the DGM and the DRM with ResNet and FCNet at their  local minimizers
	$\theta_G$ and $\theta_R$ are shown in first two rows of  Figure \ref{fig: Loss landscape for FCNet and ResNet}.
 	From the comparisons between  the left (DGM) and right (DRM) columns, 
	we can heuristically see that  the DGM has a relatively flat and smooth neighborhood while 
	the DRM seems   more rough and more oscillatory near $\theta_R$. This difference remains true  both for the fully-connected network and the ResNet.
	We change the  set of optimal parameters to
	the second set $\wt{\theta}_G$ and $\wt{\theta}_R$
	in the subfigure (c) 
	and we still see the similar observation.
	Therefore, the  visualization results we obtained  here from random directions qualitatively
	confirms our conjecture    that the DRM has  more rough landscapes near its local minimizers, while the landscapes of the DGM at local minimizers
	are relatively less rough.
 	
	\subsubsection{Understanding  difference of the roughness index for two models}
	
	\begin{figure}[htpb]
		\centering
		\includegraphics[width=2.5in]{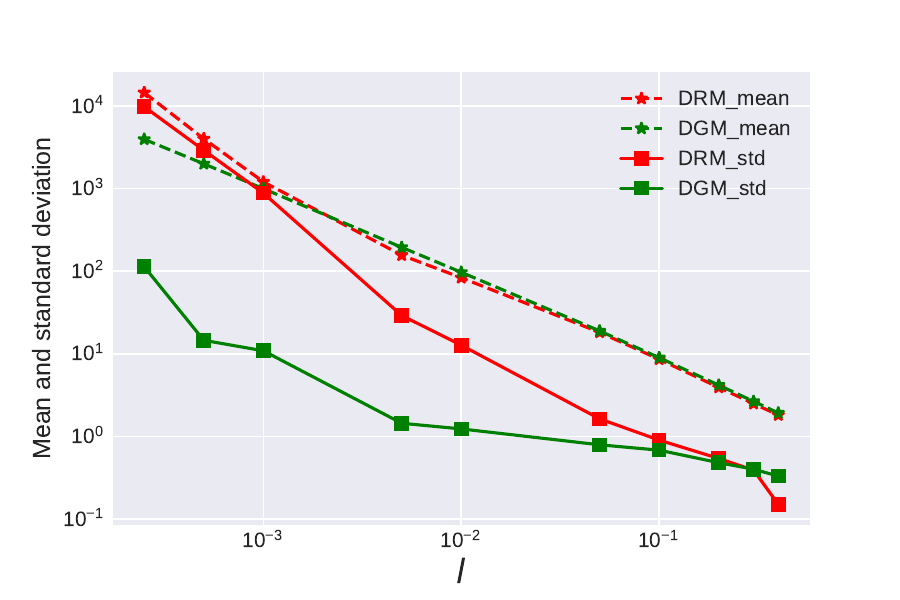}		
		\caption{The mean   $\mu$ and the std $\sigma$ in the  roughness index $\mathcal{I}=\sigma/\mu$   at $\theta_G$ and $\theta_R$ respectively, for various interval lengths $l$. $M = 100$, $m = 100$ and ResNet.
		}
		\label{fig: T for li}
	\end{figure}
	
	Recall the definition of roughness index, $\mathcal{I}=\sigma/\mu$, is the ratio of the standard deviation and the expectation of the 
	(1D) normalized TV \eqref{TT} when the loss function is projected on $M$ random directions. 
	After establishing that $\mathcal{I}$ are indeed different for the DGM and the DRM at the local minimizers $\theta_G$ and $\theta_R$, respectively. 
	we want to further check whether   the reason is from the standard deviation $\sigma$ or the expectation $\mu$. 
	Figure \ref{fig: T for li}  discovers that the difference comes  from the standard deviation $\sigma$, not the mean $\mu$.
	In fact, the means of the normalized TV across different directions are almost identical in the two models. 
	This figure strongly indicates the importance of taking account of random effect of the directions. A larger $\sigma$ means a higher
	anisotropy of the loss function in the high dimension. Therefore, we can say the higher roughness of the DRM comes from the
	more anisotropic loss function. 
	\begin{remark}
		Note that the ``anisotropy'' here has nothing to do with the eigenvalues of the Hessian matrix. Some conventional literatures use the 
		ratio of eigenvalues to represent the anisotropy for a quadratic function.  However, we have known that the roughness index is null
		for   quadratic functions. The  ``anisotropy''  refers to the uncertainty of the TV norms (the ``1D'' roughness) 
		across different directions in a high dimensional space. 
		\end{remark}

	\subsubsection{Roughness index on gradient descent path}
	
	So far we have focused on the roughness index around the local minimizer $\theta^*$ 
	(which is chosen as $\theta_G$,  $\theta_R$, $\wt{\theta}_{G}$,  $\wt{\theta}_R$ respectively) and we have well established  the distinctions 
	between the DGM and the DRM.  
	One natural question to follow  is whether this  significant distinction of roughness indices  at the local minimizers
	remain true  everywhere for the two loss functions.  The answer is no: the disparity of the roughness only appears near the local minimizers. 
	We provide the evidences in the following.
	Firstly,  we compute the RI for arbitrarily  points in the parameter space by following the standard  strategies 
	such as  Xaiver  initializations \cite{glorot2010understanding}, and two other  random samples.  Table \ref{tab: tvd for initialization} shows that the difference in the index 
	is very marginal. In fact,  we observed from this table that  the expectation  $\mu$ is nearly  $1/2l$ for almost  every direction.
This means the  1D projected  loss function is monotonic in all directions at all initial points:  the loss landscape  is essentially  non-oscillatory  almost everywhere
for random locations.
	
	\begin{table}[htbp]
		\centering
		\caption{The expectation $\mu$ and the std $\sigma$ in RI  at  random   points from  different initialisation strategies.
		Network  width $w = 4$, $l = 0.01$, $M = 100$, and $m = 20$.}\label{tab: tvd for initialization}
		\begin{tabular}{c|cc|cc}
			\toprule[1pt]
			
			\noalign{\smallskip}
			\multirow{2}*{Initialization} 
			&\multicolumn{2}{c|}{{DGM}} &\multicolumn{2}{c}{{DRM}}\\			
			&\multicolumn{1}{c}{$\mu$} & \multicolumn{1}{c|}{$\sigma$} & \multicolumn{1}{c}{$\mu$} & \multicolumn{1}{c}{$\sigma$}\\
			\noalign{\smallskip}
			\midrule[1pt]
			\noalign{\smallskip}
			\multirow{1}*{Xavier}
			& $50.00$ & $7.141\text{e-15}$ & $50.00$ & $7.105\text{e-15}$\\
			\multirow{1}*{Uniform(-1, 1)} 
			& $50.13$ & $1.340$ & $50.04$ & $0.3669$\\
			\multirow{1}*{Normal(0, 1)} 
			& $50.00$ & $1.596\text{e-15}$ & $50.10$ & $0.9939$\\
			\noalign{\smallskip}
			\bottomrule[1pt]
		\end{tabular}

	\end{table}

\begin{figure}[htpb!]
		\centering
		\includegraphics[width=2.5in]{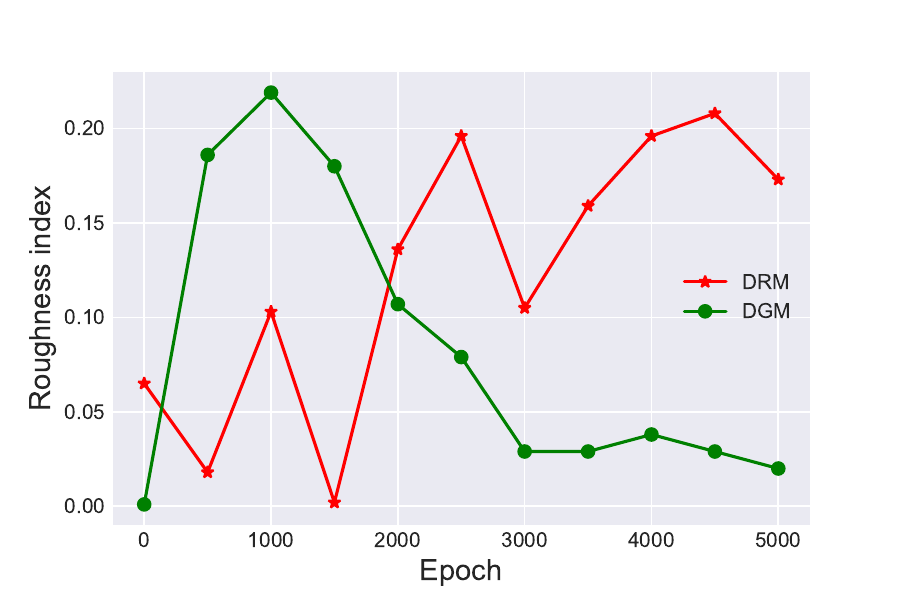}		
		\caption{The  roughness indices along  the path generated from the training process with $l  = 0.01$,  $M = 50$, $m = 10$, and the ResNet.}
		\label{fig: T_process}
	\end{figure}

	The second   evidence is from the examination of the RI along a path from an initial point to the local minimizer. 
	We first generate and save  a (gradient-descent) path obtained from the training process, then compute the roughness index 
	at a few representative points which are ordered  by the epoch.  Figure \ref{fig: T_process} presents these two 
	curves of the indices for the two models and  suggests that there is a cross-over of the roughness around at the epoch $2000$.
	Recall in Figure \ref{fig: Loss and error for widenn} which records the training process,  the training processes
	in general have already approached a vicinity of the minimizer around epoch $2000$ and after that the training is 
	to mainly  improve the accuracy further within this vicinity. By dividing the training process into these two stages,
	Figure \ref{fig: T_process} essentially tells us that in these two stages, the regions that the 
	trajectories  are exploring can be very different in terms of the roughness index.
	
	\medskip

	In summary,  by intensively examining the landscapes of   the DGM and the DRM used for  the 1D Poisson equation whose solution is smooth, 
	we provide the empirical evidences to conclude that  the DGM has a  less more rough landscape than the DRM near local minimizers  in the sense of the roughness index $\mathcal{I}$ we defined before.
	This difference could heuristically explain the reason why the DGM in general can achieve a better accuracy than the DRM,  but we have to admit 
	that a rigorous mathematical connection is still lacking here
	due to the challenge of  non-convexity.

	% In the next, we present some additional results for the cases where the solution  is not sufficient smooth and 
	% the high dimensional Poisson equation with smooth solutions.  These results are in alignment with the above conclusion.

\subsection{3D equation with a low-regularity solution}
 
  To further check  our conclusion, we consider a problem with a low-regularity solution over $\Omega = \{x \in \mathbb{R}^3: |x| < 1\}$
	\begin{equation*}
	\begin{cases}
	- \Delta u = f(x), & \; \text{in} \; \Omega,\\
	u(x) = 0, & \; \text{on} \; \partial \Omega.
	\end{cases}
	\end{equation*}
	The exact solution $u(x) = \sin \left(\frac{\pi}{2} (1 - |x|)\right)$ is continuous but not differential at the origin. Then 
$	f(x) = \frac{\pi^2}{4} \sin \left(\frac{\pi}{2} (1 - |x|)\right) + \frac{\pi}{|x|} \cos \left(\frac{\pi}{2} (1 - |x|)\right).
$
 	The solution is parametrized $u(x;\theta) = (|x| - 1) \cdot \NN(x;\theta).$
 	The ResNet is used with three residual blocks and neural width $w = 8$, thus the total number of parameters is $617$. The number of epochs is $5000$ and the batch size $N$ is $1000$. Roughness indices at the same point $\theta_G$   are recorded in Table \ref{tab: Total variation compara}. These results point to the same conclusion we had before.
  %that the roughness index  at the minimizer in the DGM is smaller than that in the DRM.
	%
	%\begin{table}[htp]
	%	\centering
	%	\begin{tabular}{|c|c|c|}
	%		\hline
	%		$M$ & $\mathcal{I} _{DGM}$ & $\mathcal{I} _{DRM}$ \\
	%		\hline
	%		50 & $0.0773$ & $0.1726$ \\
	%		\hline
	%		100 & $0.0823$ & $0.1431$ \\
	%		\hline
	%		150 & $0.0799$ & $0.1173$ \\
	%		\hline
	%		200 & $0.0851$ & $0.1175$ \\
	%		\hline
	%	\end{tabular}
	%	\caption{Roughness index of attractor for different number of random directions with $l_i = 0.1, m = 10$.}\label{tab: Total variation compara_M_3d}
	%\end{table}
	%
	%\begin{table}[htp]
	%	\centering
	%	\begin{tabular}{|c|c|c|c|}
	%		\hline
	%		$l_i$ & $m$  & $\mathcal{I}_{DGM}$ & $\mathcal{I}_{DRM}$ \\
	%		\hline
	%		0.1 &  10  & $0.0823$ & $0.1431$ \\
	%		\hline
	%		0.2 &  20  & $0.1001$ & $0.1283$ \\
	%		\hline
	%		0.3 &  30  & $0.1228$ & $0.1159$ \\
	%		\hline
	%		0.4 &  40  & $0.1442$ & $0.1516$ \\
	%		\hline
	%		0.5 &  50  & $0.1499$ & $0.1646$ \\
	%		\hline
	%	\end{tabular}
	%	\caption{Roughness index of attractor for different internal length and $m$ with $M = 100$.}\label{tab: Total variation compara_l_i_m_3d}
	%\end{table}
	%
	%\begin{table}[htp]
	%	\centering
	%	\begin{tabular}{|c|c|c|}
	%		\hline
	%		$m$ & $\mathcal{I}_{DGM}$ & $\mathcal{I}_{DRM}$ \\
	%		\hline
	%		5 & $0.0579$ & $0.1426$ \\
	%		\hline
	%		10 & $0.0823$ & $0.1431$ \\
	%		\hline
	%		15 & $0.0636$ & $0.1750$  \\ 
	%		\hline
	%		20 & $0.0675$ & $0.1283$  \\ 
	%		\hline
	%	\end{tabular}
	%	\caption{Roughness index of attractor for different grid size with $M = 100$ and $l_i = 0.1$.}\label{tab: Total variation compara_m_3d}
	%\end{table}
	
	\begin{table}[h]
			\caption{Roughness indices of $\JG$ and $\JR$ at the same point $\theta_G$ in terms of the number of random directions $M$, interval of interest $l$, and the number of grid points
		$m$ in the 3D case.}\label{tab: Total variation compara}
		\begin{center}
			\begin{tabular}{c|c|c|c|c}
				\toprule[1pt]
				$l$& $m$ & $M$ & $\mathcal{I}_{DGM}$ & $\mathcal{I}_{DRM}$  \\
				\midrule[1pt]
				\multirow{7}{*}{0.1} & 5 & 100 & 0.0579 & 0.1426 \\
				\cline{2-5}
				& \multirow{4}{*}{10} & 50 & 0.0773 & 0.1726 \\
				\cline{3-5}
				&   &  100 & 0.0823 & 0.1431 \\
				\cline{3-5}
				&   &  150 & 0.0799 & 0.1173
				\\
				\cline{3-5}
				&   &  200 & 0.0851 & 0.1175 \\
				\cline{2-5}
				& 15  &  100 & 0.0636 & 0.1750 \\
				\cline{2-5}
				& 20  &  100 & 0.0675 & 0.1283 \\
				\hline
				0.2 & 20  &  100 & 0.1001 & 0.1283 \\
				\hline
				0.3 & 30  &  100 & 0.1228 & 0.1159 \\
				\hline
				0.4 & 40  &  100 & 0.1442 & 0.1516 \\
%				\hline
%				0.5 & 50  &  100 & 0.1499 & 0.1646\\
				\bottomrule[1pt]
			\end{tabular}
		\end{center}

	\end{table}
	
	\subsection{High dimensional Poisson equation}
	Our next example is the equation \eqref{eqn: poisson equation} when $d = 10$. The ResNet is used with three residual blocks and neural width $w = 20$, thus the total number of parameters is $3601$. The number of epochs is $50000$ and the batch size $N$ is $100000$.   The relative errors in both DGM and DRM are around $\text{1e-3}$ with $50000$ epochs.
	Roughness indices of attractor in terms of the number of random directions $M$, interval of interest, and the number of grid points, are recorded in Table \ref{tab: Total variation compara high}. Again,   we observe that the roughness index   in the DGM is slightly smaller than that in the DRM.

	\begin{table}[h]
			\caption{  Roughness indices of $\JG$ and $\JR$ at the same point $\theta_G$ in terms of the number of random directions $M$, interval of interest, and the number of grid points in the 10D case.}\label{tab: Total variation compara high}
		\begin{center}
			\begin{tabular}{c|c|c|c|c}
				\toprule[1pt]
				$l$& $m$ & $M$ & $\mathcal{I}_{DGM}$ & $\mathcal{I}_{DRM}$  \\
				\midrule[1pt]
				0.025 & 20  &  100 & 0.1162 & 0.1500 \\
				\hline
				\multirow{2}{*}{0.05} & 20  &  100 & 0.0770 & 0.2126\\
				\cline{2-5}
				& 40  &  100 & 0.1015 & 0.1888\\
				\hline
				\multirow{4}{*}{0.1} & 10 & 100 & 0.1189 & 0.1420 \\
				\cline{2-5}
				& 15 & 100 & 0.1045 & 0.1497 \\
				\cline{2-5}
				& \multirow{3}{*}{20} & 50 & 0.1292 & 0.1384  \\
				\cline{3-5}
				&   &  100 & 0.1151 & 0.1763 \\
				\cline{3-5}    
				&   &  150 & 0.1092 & 0.1824 \\
				\cline{2-5}
				& 40 & 100 & 0.1124 & 0.1780 \\
				\hline
				0.2 & 20  &  100 &  0.1615 & 0.1750 \\
				\bottomrule[1pt]
			\end{tabular}
		\end{center}

	\end{table}
	
		\subsection{1D wave equation}
\begin{figure}[htpb!]
	\centering
	\includegraphics[width=2.5in]{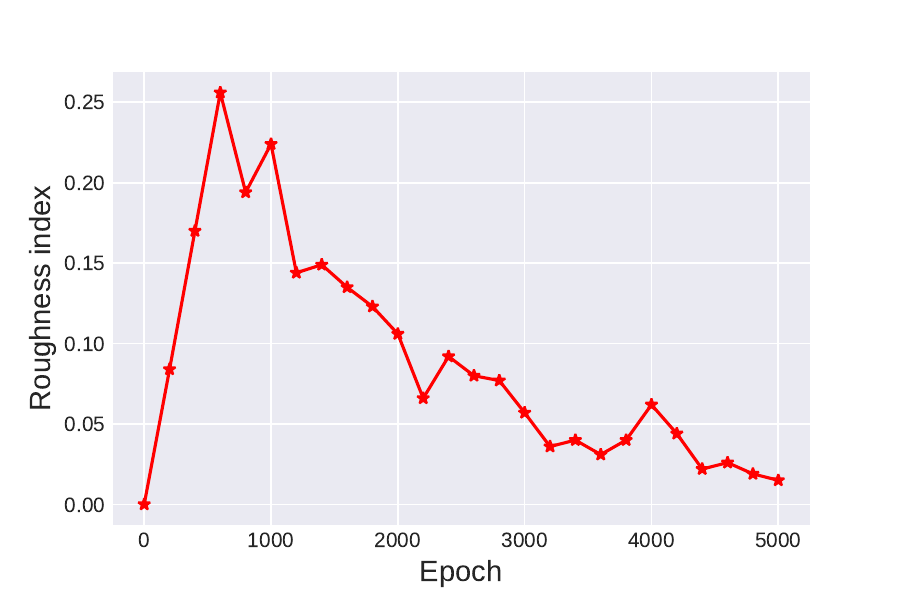}		
	\caption{The  roughness indices along  the path generated from the training process with $l  = 0.01$,  $M = 100$, $m = 40$ for 1D wave equation. The ResNet is used with one residual blocks and neural width $w = 8$.}
	\label{fig: T_process_wave}
\end{figure}

	The last example is the wave equation in one dimension:
	\begin{equation*}\label{eqn:wave}
	\begin{cases}
	u_{tt} - \Delta u = f(x), & \;  \; t\in [0, T], x\in (0,1) ,\\
	u(t, x) = 0, & \;   \; t\in  [0, T], x=0, 1\\
	u(0, x) = u_{t}(0,x)=0, & \;    x\in (0,1) 	\end{cases}
	\end{equation*}
	with the exact solution $u(x,t) = t^2   \sin(\pi x)$.
	 Similarly, the solution is parametrized by the DNN approximation
$	u(x;\theta) = t^2 x(1 - x) \cdot NN(x;\theta).
$
The Deep Ritz method   is not  applicable here because the wave equation  has no variational formulation. 
So  instead of comparing the landscapes of the  DGM and the DRM, we explore  the change of RI along
a path from a gradient  descent in training the loss function.
Figure~\ref{fig: T_process_wave} presents this  curves of the RI along with the value of the DGM loss.
We find that  the RI value along  the  path  is quite similar to that for the DGM   in Fig. \ref{fig: T_process} for the Poisson equation:
the   gradient  descent  trajectory  first go through a high RI  region  and  then  gradually decreases together with the loss.
Since the box size $l=0.01$ is used here, we can  say the gradients near the minimzer $\theta^{*}$ are all close to zero in  the neighborhood with size $l$.

\section{Concluding remarks}\label{sec:conclusion}
	In this work, we introduce a roughness index to characterize the roughness of loss function near any reference point. Through numerous experiments, we show that this quantity is particularly 
	useful for the high dimensional parameter space and  can effectively characterize the   ``roughness''   difference  between  two neural network landscapes arising from DGM and DRM.
	Our roughness index is based on the 1D normalized total variation in any specified region,  rather than the Hessian matrix at the local minimizer as a local quadratic approximation,
	so this index can be applied to both convex and non-convex landscapes.  Furthermore, 
	we propose an efficient algorithm to  compute this roughness index by randomly sampling the projection directions.
	
	In the comparison between DGM and DRM, we see significant smaller values of the roughness index for the DGM than for the DRM 
	at various local minimizers when ResNet is used.   We also discover that this difference of the roughness  
	mainly comes from the standard deviation of the directional randomness.
	By examining the roughness index along the optimization trajectory,
	we have the empirical observations  that although both are initialized 
	in a smooth region with low  RI, the RI in the DRM 
	gradually increases while the DGM has the ability to pass through a high RI region and then settle to a low RI  basin of the minimizer.
	We conjecture that this empirical observation of RI differences in the landscape  may be the reason of the
	performance  differences of using these two models in practice to solve high dimensional PDEs,
	such as the difference in the  accuracy of the numerical solution and the difficulties  of training the models.
	The last comment is  	although we propose the  roughness index
	and demonstrate its power in the background  
	of solving PDE problems, we think this roughness concept and our method of RI are also  important     in studying highly 
	non-convex  landscapes  for general machine-learning tasks.   
Particularly,   the signature pattern of {\it increasing-then-decreasing}  RI  on the  optimization path in the DGM, as shown in Figure \ref{fig: T_process} and Figure \ref{fig: T_process_wave},
 implies that  by following the gradient descent, the trajectory experiences the  ``{\it  flat-rough-flat}'' transition  when travelling the landscape.
	We conjecture that  this could  be also valid in many machine-learning tasks such as image classification problems, but the careful   empirical validations with heuristic or rigorous analysis are still
	yet under our investigation.
	
 \medskip
 \noindent \textbf{Acknowledgment.} 
The work of Chen is partially supported by National Key R\&D Program of China (No. 2022YFA1005200 and No. 2022YFA1005203), NSFC Major Research Plan -  Interpretable and General-purpose Next-generation Artificial Intelligence (No. 92270001 and No. 92270205), Anhui Center for Applied Mathematics, and the Major Project of Science \& Technology of Anhui Province (No. 202203a05020050).
This work of Du is partially supported by National Natural Science Foundation of China via grant 12271360.
The work of Zhou is partially supported by Hong Kong RGC GRF 11307319, 11308121, 11318522, and the NSFC/RGC Joint Research Scheme [RGC Project No. N-CityU102/20 and NSFC Project No. 12061160462].

% \medskip

% \section*{Acknowledgment}

%\noindent \textbf{Acknowledgment.} This work was supported by National Natural Science Foundation of China grants 11971021 (J. Chen) and 11501399 (R. Du), and Hong Kong RGC GRF grants 11305318 (X. Zhou). 

\bibliographystyle{IEEEtran}      
\bibliography{refs}

\end{document}